\documentclass[EJP,preprint]{ejpecp} 

\usepackage{enumitem}

\SHORTTITLE{Exact Asymptotics for the Harmonic Descent Chain}

\TITLE{Exact Asymptotics for the Harmonic Descent Chain and Extended Convergence Bounds for Recursive Statistics
} 

\AUTHORS{%
  Sean P.~Svihla\footnote{University of Colorado Boulder, United States of America.
    \EMAIL{sean.svihla@colorado.edu}}\orcid{0000-0002-3496-1147} 
}


\AMSSUBJ{60J80; 60J10} 
\AMSSUBJSECONDARY{60C05; 05A16; 05C05} 

\SUBMITTED{January 2, 2013} 
\ACCEPTED{December 13, 2014} 

\VOLUME{0}
\YEAR{2023}
\PAPERNUM{0}
\DOI{10.1214/YY-TN}

\ABSTRACT{
  We give an exact asymptotic expansion of the occupation probabilities $a(n,k)$ of the harmonic descent chain, introduced by D.\ Aldous, as $n\rightarrow\infty$ with $k\in\lilo{\sqrt{n}}$.
  In addition, we extend the recent results of Brandenberger et al.\ by giving a bound on the rate of convergence of sequences defined recursively by $x_n = \frac{1}{H_{n-1}}\sum_{i=1}^{n-1}\frac{x_i}{n-i}$ for more general initial conditions.
  This result enables the application of an existing central limit theorem, also by Brandenberger et al., to statistics of the critical beta-splitting tree which satisfy our more general initial conditions.
}

\DeclareMathOperator{\var}{Var}
\DeclareMathOperator{\cov}{Cov}
\DeclareMathOperator{\res}{Res}

\DeclareRobustCommand{\abs}[1]{\left|#1\right|}
\DeclareRobustCommand{\one}[0]{\mathds{1}}                
\DeclareRobustCommand{\Var}[1]{\var\left(#1\right)}       
\DeclareRobustCommand{\e}[0]{\mathbb{E}}                  
\DeclareRobustCommand{\E}[1]{\e\left(#1\right)}
\DeclareRobustCommand{\bigO}[1]{\mathrm{O}\left(#1\right)}
\DeclareRobustCommand{\lilo}[1]{\mathrm{o}\left(#1\right)}
\DeclareRobustCommand{\Res}[2]{\res_{#1}\left[#2\right]}
\DeclareRobustCommand{\Exp}[1]{\exp\left\{#1\right\}}

\newcommand{\given}{\;|\;}

\newcommand{\normal}{\calN}                         
\newcommand{\st}{_{\mathrm{st}}}                    
\newcommand{\indset}[1]{\one\!\left(#1\right)}      

\newcommand{\leaf}[1]{\left[#1\right]}
\newcommand{\parent}[1]{p\left(#1\right)}
\newcommand{\subtree}[2]{#1(#2)}

\newcommand{\rmb}{\mathrm{b}}

\newcommand{\rmX}{\mathrm{X}}
\newcommand{\rmT}{\mathrm{T}}
\newcommand{\rmN}{\mathrm{N}}
\newcommand{\rmI}{\mathrm{I}}

\newcommand{\calB}{\mathcal{B}}
\newcommand{\calN}{\mathcal{N}}

\newcommand{\calC}{\mathcal{C}}
\newcommand{\calR}{\mathcal{R}}

\renewcommand{\d}{\,d}
\newcommand{\dt}{\d t}
\newcommand{\du}{\d u}
\newcommand{\ds}{\d s}
\newcommand{\dy}{\d y}

\newcommand{\R}{\mathbb{R}}

\newenvironment{case}[2]
  {%
    \renewcommand{\proofname}{\itshape Case #1\upshape{} (#2)}
    \begin{proof}%
  }
  {\end{proof}}

\renewcommand{\i}{\mathrm{i}}                       
\renewcommand{\Im}[1]{\mathrm{Im}\left(#1\right)}   
\renewcommand{\Re}[1]{\mathrm{Re}\left(#1\right)}   

\begin{document}



\section{Introduction}

The critical beta-splitting tree is a particular case of a one-parameter family of random binary trees called beta-splitting models.
The model was introduced in~\cite{Aldous_1996} as a more realistic model of real-world phylogenetic trees, and it has seen renewed interest in recent years~\cite{Aldous_Pittel_2025,Aldous_Janson_2024a,Aldous_Janson_2024b,Aldous_Janson_2025,Aldous_Janson_Li_2024}.
Under a beta-splitting model, a random tree is generated by recursively partitioning $\{1, \dots, n\}$ according to a splitting probability distribution $q(n,i)$, $1\leq i < n$, which gives the probability that a clade of size $n$ is split into clades of size $i$ and $n-i$.
For the critical beta-splitting tree, this probability is 
\begin{equation}\label{eqn:qni-def}
    q(n,i) := \frac{n}{2H_{n-1}}\frac{1}{i(n-i)},
\end{equation}
where $H_{n-1}$ denotes the $(n-1)$th harmonic number.

The harmonic descent chain, described in~\cite{Aldous_Pittel_2025,Aldous_Janson_2024a,Aldous_Janson_2024b,Aldous_Janson_Li_2024}, is a discrete-time Markov chain on the states $\{1, 2, 3, \dots\}$ with transition probabilities 
\[
    p(j,i) = \frac{1}{(j-i)H_{j-1}}, \quad 1\leq i < j,\, j\geq 2,
\]
and $p(1,1) = 1$.
Started at state $n$, this Markov chain describes the sequence of clade sizes in the critical beta-splitting tree with $n$ leaves along a path from the root to a uniformly random leaf. 
Of interest are the ``occupation probabilities'', $a(n,k)$, defined as the probability that the chain started at state $n$ is ever in state $k$.
These characterize the asymptotic fringe distribution of the tree (see~\cite{Aldous_1991} in general), that is, the $n\rightarrow\infty$ weak limit of the tree relative to a uniformly random leaf.
In particular, it has been shown~\cite{Aldous_Janson_2025,Aldous_Janson_2024b,Aldous_Janson_Li_2024,Iksanov_2025} that the limits $\lim_{n\rightarrow\infty}a(n,k) =: a(k)$ exist, 
\begin{equation}
    a(k) = \frac{H_{k-1}}{\zeta(2)(k-1)}, \quad k\geq 2,
\end{equation}
and $a(1) = 1$.
In turn, the ``reverse harmonic descent chain'' is a Markov chain with ``upward'' transition probabilities
\[
    q^{\uparrow}(i,j) = \frac{a(j)}{a(i)}p(j,i).
\]
Letting, at each upward step $i\rightarrow j$, the sibling clade be distributed according to the critical beta-splitting distribution with $j-i$ leaves, the resulting infinite tree described by the upward chain is the extended fringe~\cite{Aldous_1991} relative to a typical leaf.

Many quantities of interest obey the recursion 
\begin{equation}\label{eqn:xn-recursion-intro}
    x_n = \frac{1}{H_{n-1}}\sum_{i=1}^{n-1}\frac{x_i}{n-i}, \quad n > M,
\end{equation}
with some set of initial conditions.
For example, if $\rmN_{n}(k)$ is the number of $k$-leaf subtrees of a critical beta-splitting tree with $n$ leaves and $\rmI_n$ is the size of the split at the root, then the conditional distribution of $\rmN_{n}(k)$ given $\rmI_n = i$ is the same as $\rmN_{i}(k) + \rmN_{n-i}(k) + \indset{n=k}$.
So, the law of total probability yields 
\[
    \frac{\E{\rmN_{n}(k)}}{n} = \frac{1}{H_{n-1}}\sum_{i=1}^{n-1}\frac{\E{\rmN_{i}(k)}/i}{n-i}
\]
for $n > k$, that is, the quantity $\E{\rmN_{n}(k)}/n = a(n,k) / k$ obeys the recursion described.

Recently, Brandenberger et al.~\cite{Brandenberger_Chin_Mossel_2026} analyzed such recursions and found bounds on the rate of convergence to the limit $x:=\lim_{n\rightarrow\infty}x_n$ when $x_n$ has the initial conditions $x_1 = \cdots = x_{k-1} = 0$ and $x_k > 0$.
Because $a(n,k) = k\E{\rmN_{n}(k)}/n$, this analysis yielded, in turn, a description of the asymptotic behavior of $a(n,k)$;
in particular, they showed~\cite[eqn.\ (1.6)]{Brandenberger_Chin_Mossel_2026}
\begin{equation}
    a(n,k) - a(k) = n^{-\gamma_* + \lilo{1}}
\end{equation}
for some $1.5 < \gamma_* < 2$.
Moreover, they showed a central limit theorem based on the contraction method (e.g., see~\cite{Régnier_1989,Rösler_1991} for this method in general) for statistics of the critical beta-splitting tree:

\begin{proposition}[Proposition 1.3 of~\cite{Brandenberger_Chin_Mossel_2026}]\label{prop:brandenberger-clt}
    Let ${(\rmX_n')}_{n\geq 1}$ and ${(\rmX_n'')}_{n\geq 1}$ be independent copies of the random sequence ${(\rmX_n)}_{n\geq 1}$, independent also of ${(\rmI_n,\rmb_n)}_{n\geq 2}$, and suppose that 
    \begin{equation}
        \rmX_n\overset{d}{=}\rmX_{\rmI_n}' + \rmX_{n-\rmI_n}'' + \rmb_n, \quad n\geq 2.
    \end{equation}
    Suppose further that there exist constants $\mu,\sigma,\epsilon,C > 0$ such that $\abs{\E{\rmX_n}-\mu n}\leq C n^{1/2-\epsilon}$ and $\abs{\Var{\rmX_n}-\sigma^2 n}\leq C n^{1-\epsilon}$ and $\e(\abs{\rmb_n}^3)\leq n^{3/2-\epsilon}$.
    Then 
    \begin{equation}
        \frac{\rmX_n - \E{\rmX_n}}{\sqrt{\Var{\rmX_n}}}\overset{d}{\longrightarrow}\normal(0,1).
    \end{equation}
\end{proposition}

\noindent See~\cite{Kolesnik_2025} as well for another example of the contraction method applied specifically to the critical beta-splitting tree.
Of interest to us in particular, Brandenberger et al.~\cite{Brandenberger_Chin_Mossel_2026} demonstrated that their analysis of $x_n$ can be used to determine that ${(\rmX_n)}_{n\geq 1}$ satisfies the conditions of the CLT when ${(\E{\rmX_n}/n)}_{n\geq 1}$ obeys this recursion and initial conditions.

Our first main result is to present an exact asymptotic expansion of $a(n,k)$ in terms of $a(k)$, as $n\rightarrow\infty$ with $k\in\lilo{\sqrt{n}}$.
Our result refines the asymptotic estimate shown in~\cite{Brandenberger_Chin_Mossel_2026} and confirms the exact value of $\gamma_*$, given numerically in~\cite{Brandenberger_Chin_Mossel_2026}, as $\gamma_* = \abs{s_1} + 1 = 1.567\dots$ where $s_1$ denotes the largest negative root of $\psi(s)-\psi(1)$.
Our second result gives convergence bounds for a sequence $x_n$ obeying the recursion (\ref{eqn:xn-recursion-intro}) with a more general set of initial conditions than in~\cite[Theorem 1.1]{Brandenberger_Chin_Mossel_2026}, in particular, $x_k\geq 0$ for all $1\leq k\leq M$ and $x_k > 0$ for at least one $k$.
For appropriate random variables, this result enables Proposition~\ref{prop:brandenberger-clt} to be applied when ${(\E{\rmX_n}/n)}_{n\geq 1}$ obeys our more general initial conditions, and as a demonstration we give a CLT\ for a random variable of interest which falls outside the scope of~\cite[Theorem 1.1]{Brandenberger_Chin_Mossel_2026}.

\subsection{Notation and Terminology}

$\one_A(x)$ denotes the indicator function over the event $A$.
We may often write, for example, $\indset{x\leq N}$ for $\one_{\{y\leq N\}}(x)$.

We use standard terminology for trees and consider ordered, unlabelled binary trees.
For a tree $T$, let $\abs{T}$ denote the number of vertices in $T$, and let $\leaf{T}$ denote the number of leaves of $T$.
For a vertex $v\in T$, let $\subtree{T}{v}$ denote the subtree of $T$ rooted at $v$, and let $\parent{v}$ denote the parent vertex of $v$.

Unless otherwise stated, $\rmT_n$ is a random critical beta-splitting tree with $n$ leaves.

Write $\psi(s)$ for the digamma function, $H_{n}=\sum_{i=1}^{n}i^{-1}$ for the $n$th harmonic number, and $\zeta(s)$ for the Riemann zeta function.

We denote that $f(n)$ and $g(n)$ are asymptotically equivalent by $f(n)\sim g(n)$.
\section{Main Results}
\subsection{Asymptotics for the Harmonic Descent Chain}

Our first result is an exact asymptotic expansion of the occupation probabilities, $a(n,k)$.
It was shown in Brandenberger et al.~\cite[eqn.\ (1.6)]{Brandenberger_Chin_Mossel_2026} that for $x_n$ obeying the recursion (\ref{eqn:xn-recursion-intro}) with initial conditions $x_n=0$ for $n < k$ and $x_k > 0$, that $x:= \lim_{n\rightarrow\infty}$ exists and
\[
    x_n - x = n^{-\gamma_*+\lilo{1}}
\]
for some $1.5 < \gamma_* < 2$.
They reported that $\gamma_* = 1.567353\dots$ numerically.
Because $a(n,k) = k\E{\rmN_n(k)}/n$ and $\E{\rmN_n(k)}/n$ follows this behavior of $x_n$, the result gives an asymptotic estimate for $a(n,k) - a(k)$.

We further refine this asymptotic estimate by providing a full asymptotic expansion for $a(n,k)-a(k)$;
moreover, this expansion analytically confirms the value of $\gamma_*$ as $\abs{s_1}+1 = 1.567\dots$, where $s_1$ is the largest negative root of $\psi(s)-\psi(1)$, and gives the exact coefficient of the leading term in the expansion. 
The proof proceeds by first noticing that an expression for $a(n,k)$ given in~\cite[eqn.\ (5.4)]{Aldous_Janson_2025},
\[
    a(n,k) = \binom{n-1}{k-1}H_{k-1}\sum_{j=0}^{n-k}\binom{n-k}{j}\frac{{(-1)}^{j}}{\psi(k+j)-\psi(1)},
\]
can be viewed as a forward difference of the function 
\[
    \frac{1}{\psi(k+x)-\psi(1)},
\]
then analyzing the corresponding Rice integral by the Cauchy residue theorem.

\begin{theorem}\label{thm:a(n,k)-asym}
    As $n\rightarrow\infty$ with $k\in\lilo{\sqrt{n}}$ and $k\geq 2$,
    \begin{equation}\label{eqn:a(n,k)-asym-gamma}
        a(n,k) - a(k) \sim \sum_{i=1}^{\infty}\frac{H_{k-1}\Gamma(k + \abs{s_i})}{\Gamma(k)\psi'(s_i)}\frac{\Gamma(n)}{\Gamma(n + 1 + \abs{s_i})},
    \end{equation}
    where $0 > s_1 > s_2 > \dots$ are the negative roots of $\psi(s)-\psi(1)$.
    Moreover, for fixed $k\geq 2$, we have 
    \begin{equation}\label{eqn:a(n,k)-asym-pow}
        a(n,k) - a(k) \sim \sum_{i=1}^{\infty}\sum_{j=1}^{\infty}c_{i,j}(k)n^{-\abs{s_i}-j}
    \end{equation}
    for some computable constants $c_{i,j}(k)$, and in particular,
    \begin{equation}\label{eqn:a(n,k)-asym-bigO}
        a(n,k) - a(k) = \frac{H_{k-1}\Gamma(k+\abs{s_1})}{\Gamma(k)\psi'(s_1)}n^{-\abs{s_1}-1} + \bigO{n^{-\abs{s_1}-2}}.
    \end{equation}
\end{theorem}

\subsection{Extension of Convergence Bounds for a Recursive Sequence}

Our second result gives convergence bounds on a sequence obeying the recursion (\ref{eqn:xn-recursion-intro}) with a more general set of initial conditions.
An analogous result was shown in Brandenberger et al.\ for the initial conditions $x_n=0$ for $n < k$ and $x_k > 0$.
We extend their result to the initial conditions $x_n=0$ for $n\leq M$ except $x_k>0$ for some $k\leq M$, that is, allowing the non-zero initial condition to occur anywhere in the first $M$ terms, rather than strictly at the first non-zero index.
As a corollary, we find that the result easily extends to arbitrary non-negative initial conditions.

The following theorem is analogous to~\cite[Theorem 1.1]{Brandenberger_Chin_Mossel_2026};
in particular, our parts (i) and (ii) are analogous to their parts (i) and (ii).

\addvspace{\topsep}\noindent\textbf{Outline of the proof.}
As in Brandenberger et al.~\cite{Brandenberger_Chin_Mossel_2026}, considering the more general recursion 
\[
    x_n = \sum_{i=1}^{n-1}p(n,i)x_i, \quad n > M,
\]
for some probability mass function (p.m.f.) $p(n,i)$, one finds a recursion for the sequence of differences $x_{n+1}-x_n$.
In our case, this recursion retains a forcing term which depends on the initial conditions.
In Lemma~\ref{lem:monotonicity}, which plays an analogous role to~\cite[Proposition 1.2]{Brandenberger_Chin_Mossel_2026}, we impose an additional restriction on $p(n,i)$ which is sufficient to control the behavior of the forcing term.
In particular, define 
\[
    s(n,i) = \sum_{j=1}^{i}q(n,j) \quad\text{and}\quad \alpha_{i,n} = \frac{s(n+1,i)}{s(n+1,M)} - \frac{s(n,i)}{s(n,M)}.
\]
We impose that for all $n$, $\alpha_{i,n}/\alpha_{i,n+1} \geq \alpha_{j,n}/\alpha_{j,n+1}$ whenever $i \geq j$, for $i,j\neq M$.
We find, then, that the behavior of the sequence falls into two distinct regimes: 
the differences are non-negative for all $n > M$, or the differences are initially non-negative but finally strictly negative, changing sign no more than once.

Theorem~\ref{thm:recursion-asymptotics}\ref{item:recursion-asymptotics-existence} then follows by showing that $q(n,i)$ obeys the additional constraint.
To this end, we define $Q(n,x) := [\psi(n)-\psi(n-x)]/[\psi(n)-\psi(n-M)]$, $A(n,x) := Q(n+1,x) - Q(n,x)$ and $R(n,x) := A(n,x)/A(n+1,x)$, so that $A(n,i) = \alpha_{i,n}$ and the result follows from $\partial_x R(n,x) > 0$.
This is equivalent to showing that
\[
    \frac{A'(n,x)}{A(n,x)}
\]
is decreasing in $n$, which we prove by showing that under a change of measure $A'(n,x)/A(n,x)$ can be seen as the expectation of a decreasing function $\Psi(t)$, defined in the proof, under a probability measure $\eta_{n,x}$ for which $\eta_{n+1,x}\geq\st\eta_{n,x}$.

To show Theorem~\ref{thm:recursion-asymptotics}\ref{item:recursion-asymptotics-error-bound}, we bound the rate of convergence, as in~\cite{Brandenberger_Chin_Mossel_2026}, using that the sequence of differences changes sign at most once.

\begin{theorem}\label{thm:recursion-asymptotics}
    Let $M\geq 1$. 
    Consider ${(x_i)}_{i=1}^{\infty}$ a sequence defined by
    \begin{equation}\label{eqn:xn-recursion-cbst}
        x_n = \frac{1}{H_{n-1}}\sum_{i=1}^{n-1}\frac{x_i}{n-i}, \quad\text{for } n>M,
    \end{equation}
    with initial conditions $x_n=0$ for $n\leq M$ except that $x_k>0$ for a $k\leq M$.
    Then we have the following:
    \begin{enumerate}[label=(\roman*)] 
        \item\label{item:recursion-asymptotics-existence} Either $x_n$ is monotone increasing for all $n>M$, or there exists $N\geq M$ such that $x_n$ is monotone increasing for $M< n \leq N$ and strictly monotone decreasing for $n>N$; 
        in particular, the limit $x:=\lim_{n\rightarrow\infty}x_n$ exists;
        \item\label{item:recursion-asymptotics-error-bound} There exists a constant $N\geq M$ such that for all $n>M$, $\abs{x_n - x} \leq \frac{4}{\ln(2)}\frac{N}{n-M}x_k$.
    \end{enumerate}
\end{theorem}

\begin{corollary}\label{cor:recursion-asymptotics}
    If ${(x_i)}_{i=1}^{\infty}$ is a sequence defined as in (\ref{eqn:xn-recursion-cbst}) but with the initial conditions $x_k\geq 0$ for all $1\leq k\leq M$ and $x_k > 0$ for at least one $k$, then $x:=\lim_{n\rightarrow\infty}x_n$ exists, and there exists $N'\geq M$ such that for all $n>M$, $\abs{x_n - x} \leq \frac{4}{\ln(2)}\frac{N'}{n-M}\sum_{k=1}^{M}x_k$.
\end{corollary}

\begin{figure}
    \centering
    \includegraphics[width=\textwidth]{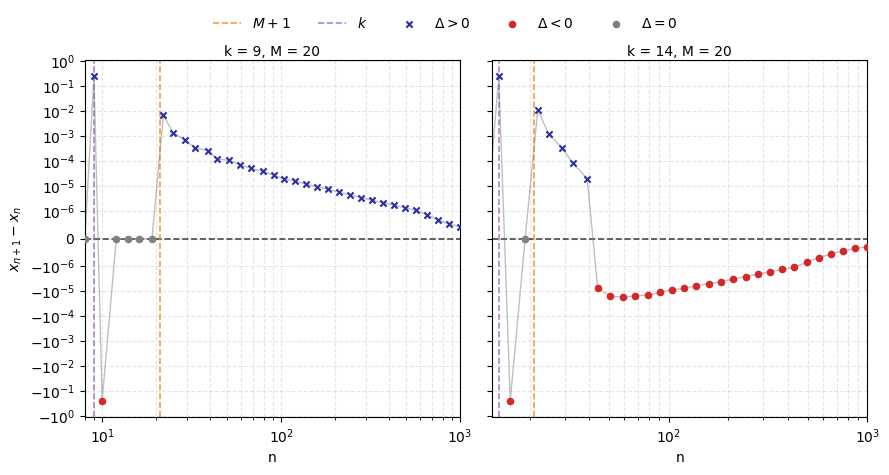}
    \caption{Illustration of tail behavior regimes.
    The left plot shows $x_{n+1}-x_n \geq 0$ for all $n>M$;
    the right plot shows $x_{n+1}-x_n$ initially non-negative, changing sign once before a strictly negative tail.
    }\label{fig:recursion}
\end{figure}

To fix ideas about Theorem~\ref{thm:recursion-asymptotics}\ref{item:recursion-asymptotics-existence}, see Figure~\ref{fig:recursion}.

\begin{remark}
    From the proof of Theorem~\ref{thm:recursion-asymptotics}\ref{item:recursion-asymptotics-error-bound}, we see that $N = M$ whenever $k = M$ or $x_n$ is monotone increasing for all $n$.
    Otherwise, $N$ is the index of the last non-negative term in $(x_{n+1}-x_n)$ and, for a particular set of initial conditions, may be found explicitly by computing the first terms of $x_n$.
\end{remark}

\subsection{Example: A Central Limit Theorem}
As an application of Theorem~\ref{thm:recursion-asymptotics}, we show a central limit theorem for a random variable of interest which does not fall in the setting of~\cite[Theorem 1.1]{Brandenberger_Chin_Mossel_2026}.

For fixed $M>0$, define $\calB_n^{(M)}$ as the minimal set of vertices $v\in\rmT_n$ which satisfies: (1) $\leaf{\subtree{\rmT_n}{v}}\leq M$ for all $v$, and (2) the collection of leaf sets of $\subtree{\rmT_n}{v}$ partitions the leaf set of $\rmT_n$.
In particular,
\begin{equation}
    \calB_n^{(M)} := \{v\given\leaf{\rmT_n(v)}\leq M < \leaf{\rmT_{n}(\parent{v})}\}
\end{equation}
Then define $\rmN_{n,k}^{(M)}$ to be the number of subtrees in $\calB_n^{(M)}$ with $k$ leaves.
Note that $\rmN_{n,k}^{(M)} = \indset{n=k}$ deterministically for $n\leq M$.

\begin{corollary}\label{cor:NnkM-CLT}
    Let $\rmX_n = \sum_{k=1}^{M}x_k\rmN_{n,k}^{(M)}$, for some constants $x_k$, be a non-constant random variable.
    Then
    \begin{equation}
        \frac{\rmX_n - \E{\rmX_n}}{\sqrt{\Var{\rmX_n}}} \overset{d}{\rightarrow} \normal(0,1).
    \end{equation}
\end{corollary}

\begin{remark}
    For example, if $x_k = k$, then $\rmX_n = n$ deterministically, and the CLT does not hold.
\end{remark}

\begin{lemma}\label{lem:|En-mu|}
    Let $E_n := n^{-1}\E{\rmX_n}$.
    There exist $\mu$ and $C$ such that $\abs{E_n - \mu} \leq Cn^{-1}$.
\end{lemma}

\begin{proof}
    Let 
    \[
        \rmX_n^+ := \sum_{k=1}^{M}x_k\indset{x_k > 0}\rmN_{n,k}^{(M)} \quad\text{and}\quad \rmX_n^- := -\sum_{k=1}^{M}x_k\indset{x_k < 0}\rmN_{n,k}^{(M)} 
    \]
    be the positive and negative parts of $\rmX_n$, respectively, so that $\rmX_n = \rmX_n^+ - \rmX_n^-$.
    Then write $E_n^+ := n^{-1}\E{\rmX_n^+}$ and $E_n^- := n^{-1}\E{\rmX_n^-}$.
    The law of total expectation yields for $n > M$
    \begin{align*}
        \E{\rmX_n^+} &= \frac{n}{2H_{n-1}}\sum_{i=1}^{n-1}\frac{\E{\rmX_i^+} + \E{\rmX_{n-i}^+}}{i(n-i)} \\
        &= \frac{n}{H_{n-1}}\sum_{i=1}^{n-1}\frac{\E{\rmX_i^+}}{i(n-i)} \\
        E_n^+ &= \frac{1}{H_{n-1}}\sum_{i=1}^{n-1}\frac{E_i^+}{n-i}
    \end{align*}
    and similarly for $\rmX_n^-$.
    $\rmX_n^+$ has the initial conditions $(x_k / k)\indset{x_k>0}$ and $\rmX_n^-$ the initial conditions $(\abs{x_k}/k)\indset{x_k<0}$ for $k=1,\dots,M$. 
    By Corollary~\ref{cor:recursion-asymptotics}, there exist $\mu^+$ and $\mu^-$ such that 
    \begin{align*}
        \abs{E_n^+ - \mu^+} &\leq \frac{4}{\ln(2)}\frac{N'}{n-M}\sum_{k=1}^{M}x_k\indset{x_k>0}\text{, and }\\
        \abs{E_n^- - \mu^-} &\leq \frac{4}{\ln(2)}\frac{N''}{n-M}\sum_{k=1}^{M}\abs{x_k}\indset{x_k<0},
    \end{align*}
    for some $N', N'' \geq M$.
    Then 
    \begin{align*}
        & \abs{E_n - (\mu^+ - \mu^-)} \leq \abs{E_n^+ - \mu^+} + \abs{E_n^- - \mu^-} \\
        & \hspace{4.5cm} \leq \frac{1}{n-M}\cdot\frac{4}{\ln(2)}\sum_{k=1}^{M}\abs{x_k}\Big[N'\indset{x_k>0} + N''\indset{x_k<0}\Big]. \qedhere
    \end{align*}
\end{proof}

\begin{lemma}\label{lem:|Vn-sigma|}
    Let $V_n := n^{-1}\Var{\rmX_n}$.
    There exist $\sigma$ and $C$ such that $|V_n - \sigma^2| \leq Cn^{-1/6}$.
\end{lemma}

\begin{proof}
    For $n > M$, the law of total variance yields 
    \begin{align*}
        V_n &= \frac{1}{n}\E{\Var{\rmX_{\rmI_n} + \rmX_{n-\rmI_n}\given\rmI_n}} + \frac{1}{n}\Var{\E{\rmX_{\rmI_n} + \rmX_{n-\rmI_n}\given\rmI_n}} \\
        &= \frac{1}{H_{n-1}}\sum_{i=1}^{n-1}\frac{V_i}{n-i} + \epsilon_n
    \end{align*}
    where 
    \[
        \epsilon_n := \frac{1}{n}\Var{\E{\rmX_{\rmI_n} + \rmX_{n-\rmI_n}\given\rmI_n}} = \bigO{\frac{1}{n}}
    \]
    by Lemma~\ref{lem:|En-mu|}.
    Since $\rmX_n$ is deterministic for $n\leq M$, we have $V_1 = \cdots = V_M = 0$.
    This recursion is identical to that of~\cite[Lemma 4.2]{Brandenberger_Chin_Mossel_2026}, and we obtain the same result.
\end{proof}

\begin{remark}
    In Lemma~\ref{lem:|Vn-sigma|},~\cite[Theorem 1.1(iii)]{Brandenberger_Chin_Mossel_2026} is already sufficient to handle the variance due to the fact that $\rmX_n$ is deterministic for $n\leq M$.
    In fact, $\rmX_n$ need not be deterministic for $n\leq M$ in order to derive a convergence bound for $V_n$, that is, $V_n$ may obey our more general initial conditions.
    For example, suppose we had found $V_n$ obeys 
    \[
        V_n = \frac{1}{H_{n-1}}\sum_{i=1}^{n-1}\frac{V_i}{n-i} + \epsilon_n
    \]
    with initial conditions $V_n \geq 0$ for $n\leq M$ and $V_k > 0$ for at least one $k\leq M$.
    Consider the related recursion 
    \begin{equation}\label{eqn:recursion-vn}
        v_n = \frac{1}{H_{n-1}}\sum_{i=1}^{n-1}\frac{v_i}{n-i},
    \end{equation}
    i.e., without the forcing term $\epsilon_n$.
    Let $V_n^{(M)}$ be the sequence obeying (\ref{eqn:recursion-vn}) for $n > M$ with initial conditions $V_n^{(M)} \geq 0$ for $n\leq M$ and $V_k > 0$ for at least one $k \leq M$.
    For each $m>M$, let $V_n^{(m)}$ obey (\ref{eqn:recursion-vn}) for $n > m$ with initial conditions $V_n^{(m)} = 0$ for $n < m$ and $V_m^{(m)} = \epsilon_m$.
    Then 
    \[
        V_n = V_n^{(M)} + \sum_{m=M+1}^{n-1}V_n^{(m)}.
    \]
    The first term, $V_n^{(M)}$, obeys $|V_n^{(M)} - \sigma_M^2| \in\bigO{{(n-M)}^{-1}}$ for some $\sigma_M^2$.
    For the rest of the terms, since each of them have the simpler initial conditions of~\cite{Brandenberger_Chin_Mossel_2026}, we may handle it using~\cite[Theorem 1.1(iii)]{Brandenberger_Chin_Mossel_2026} as in~\cite[Lemma 4.2]{Brandenberger_Chin_Mossel_2026}.
    For this reason, we did not show an analogous result to~\cite[Theorem 1.1(iii)]{Brandenberger_Chin_Mossel_2026}.
\end{remark}

\begin{proof}[Proof of Corollary~\ref{cor:NnkM-CLT}]
    Notice that $\rmN_{n,k}^{(M)} = \indset{n = k}$ for $n\leq M$.
    If we define
    \[
        \rmb_n = \sum_{k=1}^{M}x_k\big[\indset{n = k} - (\indset{\rmI_n = k} + \indset{n - \rmI_n = k})\big]\indset{2\leq n\leq M},
    \]
    we can write $\rmX_n = \rmX_{\rmI_n}' + \rmX_{n-\rmI_n}'' + \rmb_n$ for $n\geq 2$, where $\rmX_n'$ and $\rmX_n''$ are independent copies of $\rmX_n$ also independent of $(\rmI_n,\rmb_n)$.
    So, this sequence is amenable to~\cite[Proposition 1.3]{Brandenberger_Chin_Mossel_2026}.

    Using the estimates for $E_n$ and $V_n$ from Lemmas~\ref{lem:|En-mu|} and~\ref{lem:|Vn-sigma|}, the result follows from Proposition~\ref{prop:brandenberger-clt}.
\end{proof}

\section{Proof of Theorem~\ref{thm:a(n,k)-asym}}

\begin{proof}[Proof of Theorem~\ref{thm:a(n,k)-asym} (\ref{eqn:a(n,k)-asym-gamma})]
    For $k\geq 2$, we note the exact formula for $a(n,k)$ given in~\cite[eqn.\ (5.4)]{Aldous_Janson_2025}:
    \begin{equation}\label{eqn:a(n,k)-exact}
        a(n,k) = \binom{n-1}{k-1}H_{k-1}\sum_{j=0}^{n-k}\binom{n-k}{j}\frac{{(-1)}^{j}}{\psi(k+j)-\psi(1)}
    \end{equation}
    which, letting $N = n - k$, is the $N$th forward difference of
    \begin{equation}
        f_k(x) := \frac{1}{\psi(k+x)-\psi(1)}
    \end{equation}
    multiplied by $\binom{n-1}{k-1}H_{k-1}$.
    As a forward difference, we will consider its corresponding Rice integral (e.g., see~\cite{Flajolet_Sedgewick_1995}), that is, 
    \begin{equation}\label{eqn:rice-integ}
        \sum_{j=0}^{N}\binom{N}{j}{(-1)}^{j}f_k(j) = \frac{-1}{2\pi i}\oint_{\calC}\frac{\Gamma(N+1)\Gamma(-s)}{\Gamma(N+1-s)}f_k(s)\ds.
    \end{equation}
    It is required that the contour $\calC$ be positively oriented and chosen to encircle the poles located at $s = 0, 1, \dots, N$ but not the poles of $f_k(s)$, which here occur where $\psi(m + s) = \psi(1)$.

    We will first determine a contour which satisfies these criteria.
    Let $0 > s_{k,0} > s_{k,1} > \dots$ denote the roots of $\psi(k+s)-\psi(1)$ in decreasing order.
    Note $s_{k,j} = s_j - k$ for $j\geq 1$.
    It is well known that $\psi(s) - a = 0$ has only real solutions for $a\in\R$, since $\psi(s) = -\gamma + \sum_{n=0}^{\infty}\big(\frac{1}{n+1}-\frac{1}{n+z}\big)$ implies $\Im{\psi(z) - a} = \Im{z}\sum_{n=0}^{\infty}\frac{1}{\abs{n+z}^2}$.
    Then, since $\psi(s)$ is strictly increasing for $s\in\R_{>0}$, the largest pole of $f_k$ is $s_{k,0} = 1 - k$, which is negative for $k\geq 2$.
    That is, the poles of $f_k$ lie strictly on $\R_{<0}$.
    So, let $\calR$ be the rectangle with vertices $-(1/2)\pm\i$ and $(N+1/2)\pm \i$ and let $\calC = \partial\calR$ be its boundary.
    In addition, as $\psi(s)$ is strictly increasing on $(-m, -m+1)$ for $m = 1, 2, \dots$, we can be sure that $f_k$ has no poles in $(-k, 1-k)$, and that $f_k$ has exactly one pole, $s_{k,m}$ in each interval $(-m-k,-m-k+1)$ for $m\geq 1$.
    See~\cite[Chapter XII]{Whittaker_Watson_1927} or~\cite{Murty_2023} for additional details.
    In particular, the second pole, corresponding to the largest negative root of $\psi(s)-\psi(1)$, is located at $s_{k,1} = -\abs{s_1} - k = -k-0.567\dots$.

    \begin{lemma}\label{lem:zero-ordering}
        For each $m\geq 1$,
        \begin{equation}
            \abs{s_{m+1}} > \abs{s_m} + 1 + \epsilon_m,
        \end{equation}
        for some $\epsilon_m > 0$.
    \end{lemma}

    \begin{proof}
        We use the reflection formula for $\psi(s)$:
        \[
            \psi(1-s) - \psi(s) = \pi\cot(\pi s).
        \]
        As $s_m\in (-m,-m+1)$, write $s_m = -m + \delta_m$, $\delta_m\in(0,1)$.
        For $s = s_m$, the reflection formula becomes 
        \[
            \psi(1+m-\delta_m) + \gamma = \pi\cot(\pi\delta_m).
        \]
        Note that $\psi(s)$ is strictly increasing for $s > 1$, and $2 + m - \delta_{m+1} > 1 + m > 1 + m - \delta_m$.
        Since $\cot(\pi z)$ is decreasing on $(0,1)$, then 
        \[
            \cot(\pi\delta_{m+1})-\cot(\pi\delta_m) = \frac{1}{\pi}\Big[\psi(2+m-\delta_{m+1}) - \psi(1+m-\delta_m)\Big] > 0
        \]
        implies $\delta_m > \delta_{m+1}$.
        In turn 
        \[
            \abs{s_{m+1}} - \abs{s_m} = m + 1 - \delta_{m+1} - m + \delta_m = 1 + (\delta_m - \delta_{m+1}) > 1.
        \]
        We see $\epsilon_m = \delta_m - \delta_{m+1} > 0$.
    \end{proof}

    Our next step is an application of the Cauchy residue theorem.

    \begin{lemma}
        For each $m\geq 1$,
        \begin{equation}\label{eqn:deformed-integral-lem-1}
        \begin{aligned}
            \frac{-1}{2\pi\i}\oint_{\calC}K(N,s)f_k(s)\ds &= \frac{-1}{2\pi i}\int_{L_{k,m}}K(N,s)f_k(s)\ds \\
            &+ \frac{1}{\zeta(2)(k-1)}\frac{1}{\binom{N+k-1}{k-1}} + \sum_{j=1}^{\infty}\frac{\Gamma(n-k+1)\Gamma(k+\abs{s_j})}{\Gamma(n+1+\abs{s_j})\psi'(s_j)},
        \end{aligned}
        \end{equation}
        where $L_{k,m} = c_{k,m} + \i\R$ and $c_{k,m} := s_{k,m+1} + \epsilon_m$, for some $\epsilon_m > 0$ as given in Lemma~\ref{lem:zero-ordering}.
    \end{lemma}

    \begin{proof}
        We take $c_{k,m}$ and $\epsilon_m$ as defined in the statement of the lemma.
        Define $\calR_{X,Y}$ to be the rectangle with vertices $c_{k,m}\pm iY$ and $X\pm iY$ for $Y > r_0$ and $X > N+1/2$, and consider the contour $\calC_{X,Y} := \partial\calR_{X,Y}$.
        That is, we deform $\calC$ into $\calC_{X,Y}$ so that it encircles the $m+1$ largest poles of $f_k$. 
        Write
        \begin{equation}\label{eqn:K(N,s)-def}
            K(N,s) := \frac{\Gamma(N+1)\Gamma(-s)}{\Gamma(N-s+1)} = \frac{N!}{\prod_{j=0}^{N}(j-s)}
        \end{equation}
        for the kernel of the integral in (\ref{eqn:rice-integ}).
        Then 
        \begin{equation}\label{eqn:deformed-integral-finite}
        \begin{aligned}
            \frac{-1}{2\pi\i}\oint_{\calC}K(N,s)f_k(s)\ds &= \frac{-1}{2\pi\i}\oint_{\calC_{X,Y}}K(N,s)f_k(s)\ds \\
            &+ \Res{s=s_{k,0}}{K(N,s)f_k(s)} + \sum_{j=1}^{m}\Res{s=s_{k,j}}{K(N,s)f_k(s)}.
        \end{aligned}
        \end{equation}
        Denote the left segment of $\calC_{X,Y}$ from $c_{k,m}-\i Y$ to $c_{k,m} + \i Y$ by $L_{m,Y}$, and denote the remainder of the contour by $R_{X,Y}$.
        For $\abs{s} > N$, $\abs{s-j} \geq \abs{s} - N$, so 
        \begin{equation}\label{eqn:|K(N,s)|-bound}
            \abs{K(N,s)} \leq N!{(\abs{s} - N)}^{-N-1} = \bigO{\abs{s}^{-N-1}}.
        \end{equation}
        Note that $s$ is bounded away from $\R_{<0}$ and $\psi(s) = \ln(s) + \bigO{|s|^{-1}} \sim \ln\abs{s}$ as $\abs{s}\rightarrow\infty$, so that $f_k(s)\sim \ln^{-1}\abs{s}$ as $\abs{s}\rightarrow\infty$.
        As a result, the portion of (\ref{eqn:deformed-integral-finite}) along $R_{X,Y}$ converges to zero as $X,Y\rightarrow\infty$.
        Since $f_k$ has no poles and is bounded on $L_{m,Y}$, so that the integral over $L_{k,m} := \lim_{Y\rightarrow\infty}L_{m,Y}$ exists, we have 
        \begin{equation}\label{eqn:deformed-integral-infinite}
        \begin{aligned}
            \frac{-1}{2\pi\i}\oint_{\calC}K(N,s)f_k(s)\ds &= \frac{-1}{2\pi\i}\int_{L_{k,m}}K(N,s)f_k(s)\ds \\
            & + \Res{s=s_{k,0}}{K(N,s)f_k(s)} + \sum_{j=1}^{m}\Res{s=s_{k,j}}{K(N,s)f_k(s)},
        \end{aligned}
        \end{equation}
        where $L_{k,m} = c_{k,m}+\i\R$, for each $m\geq 0$.

        We handle the first residue at $s = s_{k,0}$ explicitly, since this yields the term $a(k)$.
        Since $K(N,s)$ is analytic for $\Re{s} = 1-k$, write 
        \[
            \Res{s=s_{k,0}}{K(N,s)f_k(s)} = K(N,s_{k,0})\res_{s=s_{k,0}}f_k(s)
        \]
        where
        \[
            K(N,s_{k,0}) = \frac{N!(k-2)!}{(N+k-1)!} = \frac{1}{(k-1)\binom{N+k-1}{k-1}}.
        \]
        From the Taylor series of $\psi(k+s)$ about $s_{k,0}$ and writing $u = s-s_{k,0}$,
        \[
            f_k(s) = {\Big[\psi(k+s)-\psi(1)\Big]}^{-1} = \frac{1}{\psi'(1)u}\Big[1 - \frac{\psi''(1)}{2\psi'(1)}u + \bigO{u^2}\Big].
        \]
        And, of course, $\psi'(1) = \zeta(2)$.
        That is, 
        \[
            \Res{s=s_{k,0}}{K(N,s)f_k(s)} = \frac{1}{\zeta(2)(k-1)}\frac{1}{\binom{N+k-1}{k-1}}.
        \]
        For $s = s_{k,j}$, $j\geq 1$, $K(N,s)$ is again analytic, and we have
        \[
            \Res{s=s_{k,j}}{K(N,s)f_k(s)} = K(N,s_{k,j})\Res{s=s_{k,j}}{f_k(s)} = \frac{\Gamma(N+1)\Gamma(k+\abs{s_j})}{\Gamma(N+k+1+\abs{s_j})\psi'(s_j)}.
        \]
        Substituting these values into (\ref{eqn:deformed-integral-infinite}) yields (\ref{eqn:deformed-integral-lem-1}).
    \end{proof}

    Without yet checking that this results in a valid asymptotic expansion, the above work suggests the expansion
    \begin{equation}\label{eqn:deformed-integral-infinite-with-residues}
    \begin{aligned}
        \frac{-1}{2\pi\i}\oint_{\calC}K(N,s)f_k(s)\ds &\sim \frac{1}{\zeta(2)(k-1)}\frac{1}{\binom{N+k-1}{k-1}} + \sum_{j=1}^{\infty}\frac{\Gamma(n-k+1)\Gamma(k+\abs{s_j})}{\Gamma(n+1+\abs{s_j})\psi'(s_j)}.
    \end{aligned}
    \end{equation}
    Then, substituting this expression into (\ref{eqn:a(n,k)-exact}) and writing the coefficient as $H_{k-1}\binom{n-1}{k-1} = H_{k-1}\binom{N+k-1}{k-1}$ for the first term and $H_{k-1}\binom{n-1}{k-1} = H_{k-1}\Gamma(n)/\Gamma(k)\Gamma(n-k+1)$ for the sum gives (\ref{eqn:a(n,k)-asym-gamma}).

    To check that (\ref{eqn:a(n,k)-asym-gamma}) is an asymptotic expansion, consider its the truncation at $m$ terms:
    \begin{equation}\label{eqn:deformed-integral-infinite-with-residues-truncated}
    \begin{aligned}
        a(n,k) - a(k) &= H_{k-1}\binom{n-1}{k-1}\cdot\frac{-1}{2\pi\i}\int_{L_{k,m}}K(N,s)f_k(s)\ds \\
        & \hspace{4.5cm} + \sum_{i=1}^{m}\frac{H_{k-1}\Gamma(k + \abs{s_i})}{\Gamma(k)\psi'(s_i)}\frac{\Gamma(n)}{\Gamma(n + 1 + \abs{s_i})}.
    \end{aligned}
    \end{equation}
    Note that each term in the sum is $\bigO{\ln(k)k^{\abs{s_i}}n^{-\abs{s_i}-1}}$, and for $k\in\lilo{\sqrt{n}}$, this is $\bigO{\ln(n)n^{-\abs{s_i}/2-1}}$.
    Write $s = c_{k,m} + \i y$ and recall $c_{k,m} = -\abs{s_{k,m+1}} + \epsilon_m= -\abs{s_{m+1}}-k+\epsilon_m$, for some small $\epsilon_m > 0$.
    As $\abs{y}\rightarrow\infty$, for $s\in L_{k,m}$,
    \begin{align*}
        \abs{K(N,s)} = \frac{N!}{\prod_{j=0}^{N}\abs{j-s}} &= K(N,c_{k,m})\prod_{j=0}^{N}{\Big[1 + \frac{y^2}{{(j-c_{k,m})}^2}\Big]}^{-1/2} \\
        & \leq K(N,c_{k,m}){\Big[\Big(1 + \frac{y^2}{c_{k,m}^2}\Big)\Big(1 + \frac{y^2}{{(1-c_{k,m})}^2}\Big)\Big]}^{-1/2},
    \end{align*}
    where in the final line we used that each term in the product is less than one.
    Now, notice that as $k\geq 2$, $c_{k,m} = -\abs{s_{m+1}}-k+\epsilon_m \leq -2$, and in turn, 
    $c_{k,m}^2{(1-c_{k,m})}^2 > 1$.
    It follows
    \begin{align*}
        \Big(1 + \frac{y^2}{c_{k,m}^2}\Big)\Big(1 + \frac{y^2}{{(1-c_{k,m})}^2}\Big) &= 1 + y^2\frac{c_{k,m}^2 + {(1-c_{k,m})}^2}{c_{k,m}^2{(1-c_{k,m})}^2} + y^4\frac{1}{c_{k,m}^2{(1-c_{k,m})}^2} \\
        &\geq 1 + y^2\frac{2}{c_{k,m}(1-c_{k,m})} + y^4\frac{1}{c_{k,m}^2{(1-c_{k,m})}^2} \\
        &= {\Big[1 + \frac{y^2}{c_{k,m}(1-c_{k,m})}\Big]}^2,
    \end{align*}
    so
    \[
        \abs{K(N,s)} \leq K(N,c_{k,m}){\Big[1 + \frac{y^2}{c_{k,m}(1-c_{k,m})}\Big]}^{-1}.
    \]
    Note that, as $f_k$ is bounded nearby $c_{k,m}$ on $L_{k,m}$ and $f_k(s)\sim\ln^{-1}\abs{y}$ as $\abs{y}\rightarrow\infty$, $f_k$ is bounded on $L_{k,m}$.
    So, for some $C>0$, we have 
    \begin{equation}\label{eqn:integ-bound}
        \int_{L_{k,m}}\abs{K(N,s)f_k(s)}\ds \leq K(N,c_{k,m})\int_{-\infty}^{\infty}{\Big[1 + \frac{y^2}{c_{k,m}(1-c_{k,m})}\Big]}^{-1}\dy.
    \end{equation}
    The integral here evaluates as
    \begin{equation}\label{eqn:integ-eval}
        \int_{-\infty}^{\infty}{\Big[1 + \frac{y^2}{c_{k,m}(1-c_{k,m})}\Big]}^{-1}\dy = \pi\sqrt{c_{k,m}(1-c_{k,m})} = \bigO{k}.
    \end{equation}

    \begin{lemma}\label{lem:gamma-ratio}
        As $n,k\rightarrow\infty$,
        \begin{equation}\label{eqn:coef-asym}
            H_{k-1}\binom{n-1}{k-1}\cdot K(N,c_{k,m}) = \bigO{\ln(k)k^{\abs{s_{m+1}}-\epsilon_m}n^{-\abs{s_{m+1}}-1+\epsilon_m}}.
        \end{equation}
    \end{lemma}

    \begin{proof}
        We use Stirling's approximation to handle the ratio of Gamma functions, which states $\Gamma(s) = \sqrt{2\pi} s^{s-1/2}e^{-s}[1 + \bigO{s^{-1}}]$ uniformly for $|s| \geq r_0 > 0$ and $\abs{\arg(s)} \leq \pi - \delta_0$ for fixed $r_0$ and $\delta_0 > 0$ (e.g., see~\cite[Chapter XIII]{Whittaker_Watson_1927})
        For $s = c_{k,m}$, this yields 
        \[
            \frac{\Gamma(N+1)}{\Gamma(N+1-c_{k,m})} \sim {(N+1)}^{c_{k,m}}e^{-c_{k,m}}{\Big(1 - \frac{c_{k,m}}{N+1}\Big)}^{N+1/2-c_{k,m}}.
        \]
        For large $n$, using the series expansion $-\ln(1+x)=x+\bigO{x^2}$, valid for $\abs{x}<1$, and writing $N+1/2-c_{k,m} = N+1-c_{k,m}-1/2$,
        \begin{align*}
            {\Big(1 - \frac{c_{k,m}}{N+1}\Big)}^{N+1/2-c_{k,m}} &= \Exp{-(N+1-c_{k,m}-1/2)\ln\Big(1-\frac{c_{k,m}}{N+1}\Big)} \\
            &= \Exp{(N+1-c_{k,m}-1/2)\Big(\frac{c_{k,m}}{N+1} + \bigO{\frac{c_{k,m}^2}{N^2}}\Big)} \\
            &= \Exp{c_{k,m} + \lilo{1}}.
        \end{align*}
        (Note that this step requires $k\in\lilo{\sqrt{n}}$.)
        Substituting this into the previous expression and noting ${(N+1)}^{c_{k,m}} \sim n^{c_{k,m}}$, we have 
        \begin{equation}\label{eqn:gamma-ratio-asym}
            \frac{\Gamma(N+1)}{\Gamma(N+1-c_{k,m})} \sim n^{c_{k,m}} = n^{-\abs{s_{m+1}}-k+\epsilon_m}.
        \end{equation}
        It similarly follows that
        \[
            \abs{\Gamma(-c_{k,m})} \sim \sqrt{2\pi}k^{k-1/2+\abs{s_{m+1}}-\epsilon_m}e^{-k-\abs{s_{m+1}}+\epsilon_m},
        \]
        and 
        \[
            \binom{n-1}{k-1} \sim \frac{1}{\sqrt{2\pi}}n^{k-1}e^{k-1}k^{-k+1/2}.
        \]
        Together, these estimates give (\ref{eqn:coef-asym}).
    \end{proof}

    Lastly, we combine (\ref{eqn:integ-eval}) with (\ref{eqn:coef-asym}) and (\ref{eqn:integ-bound}) to get
    \[
        H_{k-1}\binom{n-1}{k-1}\cdot\frac{-1}{2\pi\i}\int_{L_{k,m}}K(N,s)f_k(s)\ds = \bigO{\ln(k)k^{\abs{s_{m+1}}+1-\epsilon_m}n^{-\abs{s_{m+1}}-1+\epsilon_m}}.
    \]
    Since
    \[
        k^{\abs{s_{m+1}}+1-\epsilon_m}n^{-\abs{s_{m+1}}-1+\epsilon_m} = k^{\abs{s_m}}n^{-\abs{s_m}-1}\cdot\lilo{n^{(\abs{s_m}-\abs{s_{m+1}}+\epsilon_m+1)/2}},
    \]
    the result follows from $\abs{s_{m+1}} > \abs{s_m}+1+\epsilon_m$, which was shown in Lemma~\ref{lem:zero-ordering}.
\end{proof}

\begin{proof}[Proof of Theorem~\ref{thm:a(n,k)-asym} (\ref{eqn:a(n,k)-asym-pow})-(\ref{eqn:a(n,k)-asym-bigO})] 
    To obtain (\ref{eqn:a(n,k)-asym-pow}), we have the following expansion for the ratio $\Gamma(n)/\Gamma(n+1+\abs{s_i})$~\cite[eqn.\ (12)]{Erdélyi_Tricomi_1951}: 
    \[
        \frac{\Gamma(n)}{\Gamma(n+1+\abs{s_i})} \sim \sum_{j=0}^{\infty}C_{i,j}n^{-\abs{s_i}-1-j}
    \]
    where the coefficients $C_{i,j}$ are given in~\cite{Erdélyi_Tricomi_1951}.
    Substituting this into (\ref{eqn:a(n,k)-asym-gamma}) yields
    \[
        a(n,k) - a(k) \sim \sum_{i=1}^{\infty}\sum_{j=1}^{\infty}c_{i,j}(k)n^{-\abs{s_i}-j},
    \]
    where 
    \[
        c_{i,j}(k) := \frac{H_{k-1}\Gamma(k+\abs{s_i})}{\Gamma(k)\psi'(s_i)}C_{i,j-1}.
    \]
    Lastly, (\ref{eqn:a(n,k)-asym-bigO}) follows from the definition of $c_{i,j}(k)$, and in particular that $C_{i,0}=1$, and noticing that $s_1 = -0.567\dots$ and $s_2 = -1.628\dots$, so that $-\abs{s_1}-2 > -\abs{s_2}-1$ and the remainder term is $\bigO{n^{-\abs{s_1}-2}}$.
\end{proof}

\section{Proof of Theorem~\ref{thm:recursion-asymptotics} and Corollary~\ref{cor:recursion-asymptotics}}\label{sec:recursion-asymptotics}

The next result is analogous to~\cite[Proposition 1.2]{Brandenberger_Chin_Mossel_2026} and may be considered an extension of it.
The proof will be broadly similar and borrow many of its techniques.

\begin{lemma}\label{lem:monotonicity}
    Consider the sequence ${(x_i)}_{i=1}^{\infty}$ with $x_1, \dots, x_{k-1}, x_{k+1}, \dots, x_M = 0$ and $x_k > 0$ for $1\leq k\leq M$, and 
    \begin{equation}\label{eqn:xn-recursion}
        x_n = \sum_{i=1}^{n-1}p(n,i)x_i \quad\text{for } n > M,
    \end{equation}
    and define $s(n,i) := \sum_{j=1}^{i}p(n,j)$ and 
    \begin{equation}
        \alpha_{i,n} := \frac{s(n+1,i)}{s(n+1,M)} - \frac{s(n,i)}{s(n,M)}.
    \end{equation}
    Further suppose
    \begin{enumerate}[label=(\roman*)] 
        \item $p(n,i) > 0$ and $\sum_{i=1}^{n-1}p(n,i) = 1$;
        \item $p(n,i)/p(n+1,i)$ is non-decreasing in $i$ for every $n$; 
        \item for every $n$, $\alpha_{i,n}/\alpha_{i,n+1} \geq \alpha_{j,n}/\alpha_{j,n+1}$ whenever $i \geq j$, for $i,j\neq M$.
    \end{enumerate}
    Then either $x_n$ is monotone increasing for all $n > M$, or there exists $N\geq M$ such that $x_n$ is monotone increasing for $M<n\leq N$ and strictly monotone decreasing for $n > M$; in particular, the limit $x:= \lim_{n\rightarrow\infty}x_n$ exists.
\end{lemma}

\begin{proof}
    For $M\geq 2$, we will assume that $k < M$, since the case $k = M$ is simply the original result by Brandenberger et al.~\cite{Brandenberger_Chin_Mossel_2026}.

    First suppose $M,k\geq 2$ and $k < M$.
    Summation by parts yields 
    \[
        0 = x_n - x_{n-1} + \sum_{i=1}^{n-2}s(n,i)(x_{i+1}-x_i).
    \]
    Separating the $i = k-1$, $k$, and $M$ terms, we have
    \begin{align*}
        0 &= (x_n - x_{n-1})\frac{1}{s(n,M)} + \sum_{i=M+1}^{n-2}(x_{i+1} - x_i)\frac{s(n,i)}{s(n,M)} + (x_{M+1} - x_M) \\
        & \hspace{5.5cm} + (x_k - x_{k-1})\frac{s(n,k-1)}{s(n,M)} + (x_{k+1} - x_k)\frac{s(n,k)}{s(n,M)}
    \end{align*}
    for $n > M+1$ and 
    \begin{align*}
        0 &= (x_{n+1} - x_n)\frac{1}{s(n+1,M)} + \sum_{i=M+1}^{n-1}(x_{i+1} - x_i)\frac{s(n+1,i)}{s(n+1,M)} + (x_{M+1} - x_M) \\
        & \hspace{5cm} + (x_k - x_{k-1})\frac{s(n+1,k-1)}{s(n+1,M)} + (x_{k+1} - x_k)\frac{s(n+1,k)}{s(n+1,M)}
    \end{align*}
    for $n\geq M+1$.
    Subtracting the first expression from the second, we obtain 
    \begin{equation}\label{eqn:xi-difference-recursion}
        0 = (x_{n+1} - x_n)\frac{1}{s(n+1,M)} - \sum_{i=M+1}^{n-1}(x_{i+1} - x_i)\abs{\alpha_{i,n}} + x_k(\alpha_{k-1,n} - \alpha_{k,n}).
    \end{equation}
    Note $\alpha_{i,n} < 0$ for $i>M$ as shown in~\cite{Brandenberger_Chin_Mossel_2026} following from assumption (ii).
    We find by similar reasoning that $\alpha_{i,n} > 0$ for $i<M$.

    Since $x_k > 0$ and $p(n,i) > 0$ by assumption, is it easy to see from (\ref{eqn:xn-recursion}) that $x_n > 0$ for all $n>M$.
    Similarly, suppose $x_i \leq C$ for all $i < n$ and some $C > 0$.
    Then 
    \[
        x_n \leq C\sum_{i=1}^{n-1}p(n,i) = C,
    \]
    so $x_n \leq x_k$ for all $n$ by induction.
    Since the sequence is bounded above and below, we need only show that the tail of the sequence is monotone in order to conclude convergence.

    We can study the behavior of the sequence by splitting over cases depending on the sign of $x_k(\alpha_{k-1,n} - \alpha_{k,n})$ at step $n = M+1$: 

    \begin{case}{1}{$\alpha_{k-1,M+1} - \alpha_{k,M+1} = 0$}
        Note that in this case, $x_{M+2} - x_{M+1} = 0$ as well, and the recursion takes the form 
        \[
            0 = (x_{n+1} - x_n)\frac{1}{s(n+1,M)} - \sum_{i=M+2}^{n-1}(x_{i+1} - x_i)\abs{\alpha_{i,n}} + x_k(\alpha_{k-1,n} - \alpha_{k,n})
        \]
        for $n\geq M+2$.
        One can compute from (\ref{eqn:xn-recursion}) that the first non-zero term after $n = k$ is $x_{M+1} = p(M+1,k)x_k$.
        If $\alpha_{k-1,n} - \alpha_{k,n} = 0$ for all $n\geq M+2$, then by induction all the differences are zero, and $x_n = p(M+1,k)x_k$ for all $n\geq M+1$.
        Otherwise, $\alpha_{k-1,n} - \alpha_{k,n} \neq 0$ at some step $N$, at which point 
        \[
            0 = (x_{n+1} - x_n)\frac{1}{s(n+1,M)} - \sum_{i=N}^{n-1}(x_{i+1} - x_i)\abs{\alpha_{i,n}} + x_k(\alpha_{k-1,n} - \alpha_{k,n})
        \]
        for $i\geq N$, and the subsequent terms are handled under case 2, since by assumption (iii), if $\alpha_{i,n} - \alpha_{j,n} \geq 0$ for $i<j$ and $i, j\neq M$, then 
        \[
            \alpha_{i,n+1} - \alpha_{j,n+1} = \alpha_{i,n}\frac{\alpha_{i,n+1}}{\alpha_{i,n}} - \alpha_{j,n}\frac{\alpha_{j,n+1}}{\alpha_{j,n}} \geq 0.
        \]
        In particular, if $\alpha_{k-1,n} - \alpha_{k,n} \geq 0$, then $\alpha_{k-1,n+1} - \alpha_{k,n+1} \geq 0$ as well, so this first non-zero $\alpha_{k-1,n+1} - \alpha_{k,n+1}$ must, in fact, be positive, and the corresponding difference negative.
    \end{case}

    \begin{case}{2}{$\alpha_{k-1,M+1} - \alpha_{k,M+1} > 0$}
        Similar to before, by assumption (iii), if $\alpha_{i,n} - \alpha_{j,n} > 0$ for $i<j$ and $i, j\neq M$, then 
        \[
            \alpha_{i,n+1} - \alpha_{j,n+1} = \alpha_{i,n}\frac{\alpha_{i,n+1}}{\alpha_{i,n}} - \alpha_{j,n}\frac{\alpha_{j,n+1}}{\alpha_{j,n}} > 0.
        \]
        In particular, if $\alpha_{k-1,n} - \alpha_{k,n} > 0$, then $\alpha_{k-1,n+1} - \alpha_{k,n+1} > 0$ as well.
        It is clear from (\ref{eqn:xi-difference-recursion}) that $x_{M+2}-x_{M+1} < 0$ and by induction the remaining differences must be negative as well.
        The tail of the sequence is strictly monotone decreasing.
    \end{case}

    \begin{case}{3}{$\alpha_{k-1,M+1} - \alpha_{k,M+1} < 0$}
        Note that in this case, $x_{M+2}-x_{M+1} > 0$.
        If all subsequent differences are non-negative as well, then the tail of the sequence is monotone increasing.
        If not, then at some step $N$, $x_{N+1}-x_{N} < 0$, and 
        \[
            0 = (x_{N+1}-x_N)\frac{1}{s(N+1,M)} - \sum_{i=M+1}^{N-1}(x_{i+1}-x_i)\abs{\alpha_{i,N}} + x_k(\alpha_{k-1,N}-\alpha_{k,N}).
        \]
        implies $x_k(\alpha_{k-1,N} - \alpha_{k,N}) > 0$ and in turn $\alpha_{k-1,N} - \alpha_{k,N} > 0$.
    
        We note two preliminary results regarding $\alpha_{i,n}$ which follow from assumption (iii).
        First, by the induction used in the previous case, $\alpha_{k-1,n} - \alpha_{k,n} > 0$ for all $n\geq N$.
        Second, we find
        \begin{equation}\label{eqn:aux-ratio-ineq-1}
            \frac{\alpha_{k-1,n} - \alpha_{k,n}}{\alpha_{k-1,N} - \alpha_{k,N}} \geq \frac{\abs{\alpha_{i,n}}}{\abs{\alpha_{i,N}}}
        \end{equation}
        for $k<M<i$ and $n>N$.
        In particular, showing the result for $n$ and $n-1$, this is equivalent to $\abs{\alpha_{i,n-1}}(\alpha_{k-1,n} - \alpha_{k,n}) - \abs{\alpha_{i,n}}(\alpha_{k-1,n-1}-\alpha_{k,n-1}) \geq 0$, and
        \begin{align*}
            & \abs{\alpha_{i,n-1}}(\alpha_{k-1,n} - \alpha_{k,n}) - \abs{\alpha_{i,n}}(\alpha_{k-1,n-1}-\alpha_{k,n-1}) \\
            & \hspace{2cm} = \abs{\alpha_{i,n}}\alpha_{k,n-1} - \abs{\alpha_{i,n-1}}\alpha_{k,n} + \abs{\alpha_{i,n-1}}\alpha_{k-1,n} - \abs{\alpha_{i,n}}\alpha_{k-1,n-1} \\
            & \hspace{2cm} = \abs{\alpha_{i,n}}\alpha_{k,n}\Big(\frac{\alpha_{k,n-1}}{\alpha_{k,n}} - \frac{\alpha_{i,n-1}}{\alpha_{i,n}}\Big) + \abs{\alpha_{i,n}}\alpha_{k,n}\Big(\frac{\alpha_{i,n-1}}{\alpha_{i,n}} - \frac{\alpha_{k-1,n-1}}{\alpha_{k-1,n}}\Big) \\
            & \hspace{2cm} = \abs{\alpha_{i,n}}\Big[(\alpha_{k,n} - \alpha_{k-1,n})\Big(\frac{\alpha_{k,n-1}}{\alpha_{k,n}} - \frac{\alpha_{i,n-1}}{\alpha_{i,n}}\Big) + \alpha_{k-1,n}\Big(\frac{\alpha_{k,n-1}}{\alpha_{k,n}} - \frac{\alpha_{k-1,n-1}}{\alpha_{k-1,n}}\Big)\Big]
        \end{align*}
        which is non-negative by assumption (iii) and the result that $\alpha_{k-1,n}-\alpha_{k,n} > 0$ for $n\geq N$ and $\alpha_{k-1,n} > 0$.
        This will justify our inductive step, and by the same logic we obtain the base case when $n = N+1$.
        Supposing this holds up to $n-1$ for some $n$, then,
        \begin{align*}
            \frac{\alpha_{k-1,n} - \alpha_{k,n}}{\alpha_{k-1,N} - \alpha_{k,N}} - \frac{\abs{\alpha_{i,n}}}{\abs{\alpha_{i,N}}} &= \frac{\alpha_{k-1,n} - \alpha_{k,n}}{\alpha_{k-1,n-1} - \alpha_{k,n-1}}\cdot\frac{\alpha_{k-1,n-1} - \alpha_{k,n-1}}{\alpha_{k-1,n-2} - \alpha_{k,n-2}}\cdots\frac{\alpha_{k-1,N+1} - \alpha_{k,N+1}}{\alpha_{k-1,N} - \alpha_{k,N}} \\
            & \qquad - \frac{\abs{\alpha_{i,n}}}{\abs{\alpha_{i,n-1}}}\cdot\frac{\abs{\alpha_{i,n-1}}}{\abs{\alpha_{i,n-2}}}\cdots\frac{\abs{\alpha_{i,N+1}}}{\abs{\alpha_{i,N}}}
        \end{align*}
        which is non-negative, so the result holds by induction.

        Now suppose $x_{i+1}-x_i \geq 0$ for $M+1 \leq i < N$ and $<0$ for $N\leq i \leq n-1$.
        Then 
        \[
            0 = (x_{n+1} - x_n)\frac{1}{s(n+1,M)} - \sum_{i=M+1}^{N-1}(x_{i+1} - x_i)\abs{\alpha_{i,n}} - \sum_{i=N}^{n-1}(x_{i+1} - x_i)\abs{\alpha_{i,n}} + x_k(\alpha_{k-1,n} - \alpha_{k,n})
        \]
        where the sum over $i=N,\dots,n-1$ is positive, and 
        \begin{align*}
            & x_k(\alpha_{k-1,n} - \alpha_{k,n}) - \sum_{i=M+1}^{N-1}(x_{i+1}-x_i)\abs{\alpha_{i,n}} \\
            & \qquad \qquad = x_k(\alpha_{k-1,N}-\alpha_{k,N})\frac{\alpha_{k-1,n}-\alpha_{k,n}}{\alpha_{k-1,N}- \alpha_{k,N}} - \sum_{i=M+1}^{N-1}(x_{i+1}-x_i)\abs{\alpha_{i,N}}\frac{\abs{\alpha_{i,n}}}{\abs{\alpha_{i,N}}} > 0
        \end{align*}
        by (\ref{eqn:aux-ratio-ineq-1}) and the assumption that 
        \[
            -(x_{N+1}-x_N)\frac{1}{s(N+1,M)} = x_k(\alpha_{k-1,N}-\alpha_{k,N}) - \sum_{i=M+1}^{N-1}(x_{i+1}-x_i)\abs{\alpha_{i,N}} > 0.
        \]
        So, the tail of the sequence after step $N$ is strictly monotone decreasing.
    \end{case}

    Next suppose $M\geq 2$ and $k=1$.
    Because there is no $x_k - x_{k-1}$ term, we have instead the recursion 
    \[
        0 = (x_{n+1} - x_n)\frac{1}{s(n+1,M)} - \sum_{i=M+1}^{n-1}(x_{i+1}-x_i)\abs{\alpha_{i,n}} - x_1\alpha_{1,n}.
    \]
    Since $\alpha_{i,n} > 0$ for $i<M$ and all $n$ and $x_1 > 0$, we know $-x_1\alpha_{1,n} < 0$.
    Following logic similar to that of case 3 above for a negative forcing term, we have $x_{n+1} - x_n > 0$ for all $n\geq M+1$ by induction, so again the tail of the sequence is strictly monotone increasing.
    (In particular, we do not need to consider the forcing term changing sign, since $-x_1\alpha_{1,n} < 0$ for all $n\geq M+1$.)

    Lastly, suppose $k = M = 1$.
    Now the $x_{k+1}-x_k$ term in the original recursion is just the $x_{M+1} - x_M$ term cancelled out, since now $k$ and $M$ coincide, so 
    \[
        0 = (x_{n+1} - x_n)\frac{1}{s(n+1,1)} - \sum_{i=2}^{n-1}(x_{i+1}-x_i)\abs{\alpha_{i,n}}.
    \]
    Following the logic of case 1 for a zero forcing term, we see that $x_n = x_1$ for $n > 1$.
\end{proof}

\begin{proof}[Proof of Theorem~\ref{thm:recursion-asymptotics}\ref{item:recursion-asymptotics-existence}]
    It was shown in~\cite{Brandenberger_Chin_Mossel_2026} that assumption (ii) of Lemma~\ref{lem:monotonicity} is satisfied when $p(n,i) = {[H_{n-1}(n-i)]}^{-1}$.
    It will suffice to show assumption (iii) is satisfied as well.

    Define $Q(n,x) := [\psi(n)-\psi(n-x)]/[\psi(n)-\psi(n-M)]$, $A(n,x) := Q(n+1,x) - Q(n,x)$ and $R(n,x) := A(n,x)/A(n+1,x)$.
    We wish to show $R'(n,x) = \partial_x R(n,x) > 0$.
    Because $R(n,x)$ is undefined at $x = M$, we will first show $R'(n,x) > 0$ for $x < M$ and $x > M$ separately, then verify that $R(n,x)$ has a continuous extension to $x = M$.
    To begin, notice 
    \[
        R'(n,x) = \frac{A(n,x)A(n+1,x)}{{[A(n+1,x)]}^2}\Big[\frac{A'(n,x)}{A(n,x)} - \frac{A'(n+1,x)}{A(n+1,x)}\Big]
    \]
    so, since $A(n,x)$ and $A(n+1,x)$ have the same sign for all $x$, this is equivalent to showing 
    \[
        \frac{A'(n,x)}{A(n,x)} - \frac{A'(n+1,x)}{A(n+1,x)} > 0.
    \]

    We first show that this ratio can be viewed as a ratio of covariances under a probability measure $\mu_n$.

    \begin{lemma}
        For $g_x(t)$ and $\mu_n$ as defined in (\ref{eqn:g_x(t)-def}) and (\ref{eqn:mu_n-def}),
        \begin{equation}\label{eqn:A'(n,x)/A(n,x)-as-covariance}
            \frac{A'(n,x)}{A(n,x)} = \frac{\cov_{\mu_n}(\partial_x g_x(\tau),\tau)}{\cov_{\mu_n}(g_x(\tau),\tau)}.
        \end{equation}
    \end{lemma}

    \begin{proof}
        Note the integral definition of the digamma function:
        \[
            \psi(n) = -\gamma + \int_{0}^{1}\Big(\frac{1 - t^n}{1-t}\Big)\dt \quad\Longrightarrow\quad \psi(n) - \psi(n-x) = \int_{0}^{1}t^{n-1}\Big(\frac{t^{-x}-1}{1-t}\Big)\dt.
        \]
        To begin, we express the function $Q(n,x)$ as an expectation:
        Define the weight function $w_n(t) := t^{n-1}\big(\frac{t^{-M}-1}{1-t}\big)$ and its integral $W_n := \int_{0}^{1}w_n(t)\dt = \psi(n)-\psi(n-M)$.
        Then
        \[
            Q(n,x) = \int_{0}^{1}t^{n-1}\Big(\frac{t^{-x}-1}{1-t}\Big)\frac{dt}{W_n} = \int_{0}^{1}\Big(\frac{t^{-x}-1}{t^{-M}-1}\Big)\frac{w_n(t)dt}{W_n} = \e_{\mu_n}(g_x(\tau))
        \]
        where
        \begin{equation}\label{eqn:g_x(t)-def}
            g_x(t) := \frac{t^{-x}-1}{t^{-M}-1},
        \end{equation}
        and
        \begin{equation}\label{eqn:mu_n-def}
            d\mu_n := \frac{w_n(t)}{W_n}dt
        \end{equation}
        is a probability measure with density function $f_n(t) = w_n(t)/W_n$.
        Since $w_{n+1}(t) = t\,w_n(t)$, $W_{n+1}$ and $W_n$ are related by
        \[
            W_{n+1} = \int_{0}^{1}t\,w_n(t)\dt = W_n\int_{0}^{1}t\frac{w_n(t)}{W_n}\dt = W_n\e_{\mu_n}(\tau)
        \]
        and, in turn,
        \[
            d\mu_{n+1} = \frac{t\,w_n(t)}{W_{n+1}}dt = \frac{t}{\e_{\mu_n}(\tau)}d\mu_n.
        \]
        We note that $\mu_{n+1}$ is the $\tau$-size bias of $\mu_{n}$, but we do not use this fact explicitly.
        In particular,
        \[
            \e_{\mu_{n+1}}(g_x(\tau)) = \frac{\e_{\mu_{n}}(\tau g_x(\tau))}{\e_{\mu_n}(\tau)}.
        \]
        It follows, since $A(n,x) = Q(n+1,x) - Q(n,x)$, that
        \[
            A(n,x) 
            = \e_{\mu_{n+1}}(g_x(\tau)) - \e_{\mu_n}(g_x(\tau))
            = \frac{\cov_{\mu_n}(g_x(\tau),\tau)}{\e_{\mu_n}(\tau)},
        \]
        and similarly
        \[
            A'(n,x) = \frac{\cov_{\mu_n}(\partial_x g_x(\tau),\tau)}{\e_{\mu_n}(\tau)}. \qedhere
        \]
    \end{proof}

    From here the proof proceeds as follows: We express the ratio $A'(n,x)/A(n,x)$ as $\e_{\eta_{n,x}}(\Psi(t))$ ($\Psi$ and $\eta_{n,x}$ defined below), and then show that $\Psi(t)$ is decreasing and $\eta_{n+1,x}\geq\st\eta_{n,x}$ ($\eta_{n+1,x}$ stochastically dominates $\eta_{n,x}$), from which the result follows (e.g., see~\cite[Theorem 1.A.3]{Shaked_Shanthikumar_2007}).

    \begin{lemma}
        For $\eta_{n,x}$ as defined in (\ref{eqn:eta_xt-def}), $\eta_{n+1,x}\geq\st\eta_{n,x}$, and $A'(n,x)/A(n,x) = \e_{\eta_{n,x}}(\Psi(\tau))$ with $\Psi(t) = \partial_x \ln\abs{g'_x(t)}$.
    \end{lemma}

    \begin{proof}
        Define $E_n := \e_{\mu_n}(\tau)$.
        Writing $f(t) = \int_{0}^{t}f'(u)\du$ and changing the order of integration yields
        \begin{align}
            \cov_{\mu_n}(f(t),t) &= \int_{0}^{1}f(t)(t-E_n)d\mu_n \notag \\
            &= \int_{0}^{1}f'(u)\Big(\int_{0}^{u}(E_n-t)d\mu_n\Big)\du \notag \\
            &= \int_{0}^{1}f'(u)I_n(u)\du \label{eqn:cov-mu_n-change-of-measure-identity}
        \end{align}
        where 
        \begin{equation}\label{eqn:In-def1}
            I_n(u) := \int_{0}^{u}(E_n - t)d\mu_n(t).
        \end{equation}
        We have used that 
        \[
            \int_{0}^{1}(t - E_n)d\mu_n(t) = \int_{0}^{u}(t - E_n)d\mu_n(t) + \int_{u}^{1}(t - E_n)d\mu_n(t) = 0,
        \]
        which yields the following equivalent definition, which will become useful shortly:
        \begin{equation}\label{eqn:In-def2}
            I_n(u) := -\int_{u}^{1}(E_n - t)d\mu_n(t).
        \end{equation}
        Since $I_n(u) = 0$ for $u = 0, 1$, increasing on $0\leq u\leq E_n$ decreasing on the rest of the interval, $I_n(u) > 0$ for $0<u<1$.
        This transform is inspired by the method of zero-biasing found used with Stein's method.
        See, e.g.,~\cite{Goldstein_Reinert_1997} or~\cite[\S 2.3.3]{Chen_Goldstein_Shao_2011}.

        When $f(t) = g_x(t)$, we have 
        \begin{equation}\label{eqn:g'_x(t)}
            g'_x(t) = \frac{t^{M-x-1}}{{(t^M-1)}^2}\Big[x(t^{M}-1) - M(t^{x}-1)\Big].
        \end{equation}
        Of course, the leading factor is positive for $0<t<1$.
        The second factor equals zero at $t=1$, and we find its derivative with respect to $t$,
        \[
            \partial_t\Big[x(t^{M}-1) - M(t^{x}-1)\Big] = xMt^{x-1}\Big[t^{M-x}-1\Big],
        \]
        is negative if $x<M$ and positive if $x>M$, for all $0<t<1$.
        That is, the term is either monotone increasing or monotone decreasing towards zero.
        It follows that $g'_x(t)$ has the same sign for all $0<t<1$, and in particular $g'_x(t) > 0$ when $x < M$ and $g'_x(t) < 0$ when $x > M$.
        So, writing $K_{n,x} := \cov_{\mu_n}(g_x(\tau),\tau)$,
        \begin{equation}\label{eqn:eta_xt-def}
            d\eta_{n,x} := \frac{g'_x(t)I_n(t)dt}{K_{n,x}} = \frac{\abs{g'_x(t)}I_n(t)dt}{\int_{0}^{1}\abs{g'_x(u)}I_n(u)\du}
        \end{equation}
        is a valid probability measure with density function $\phi_{n,x}(t) := g'_x(t)I_n(t)/K_{n,x}$.

        To show $\eta_{n+1,x}\geq\st\eta_{n,x}$, it is sufficient to show that $\eta_{n+1,x}$ dominates $\eta_{n,x}$ in likelihood ratio~\cite[Theorem 1.C.1]{Shaked_Shanthikumar_2007}.
        Notice that the ratio of the densities of $\eta_{n+1,x}$ and $\eta_{n,x}$ is 
        \[
            \frac{\phi_{n+1,x}(t)}{\phi_{n,x}(t)} = \frac{K_{n,x}}{K_{n+1,x}}\cdot\frac{I_{n+1}(t)}{I_{n}(t)},
        \]
        so we need only show that $I_{n+1}(t)/I_n(t)$ is increasing in $t$.
        The sign of the derivative is the same as the sign of 
        \begin{align*}
            I'_{n+1}(t)I_n(t) - I_{n+1}(t)I'_n(t) &= \Big[(E_{n+1} - t)\frac{t f_n(t)}{E_n}\Big]I_n(t) - I_{n+1}(t)\Big[(E_n - t)f_n(t)\Big] \\
            &= \frac{f_n(t)}{E_n}\Big[t(E_{n+1} - t)I_n(t) - E_n(E_n - t)I_{n+1}(t)\Big].
        \end{align*}
        Of course, the leading factor is positive, since $f_n$ is a probability density and $E_n$ is the mean of a random variable that is positive almost surely.
        Define $H(t) := t(E_{n+1}-t)I_n(t) - E_n(E_n-t)I_{n+1}$.
        We show in the following cases that $H(t) > 0$, from which $\eta_{n+1,x}\geq\st\eta_{n,x}$:

        \begin{case}{1}{$0<t\leq E_n$}
            Write $H(t)$ under a single integral using the first definition of $I_n(t)$ in (\ref{eqn:In-def1}):
            \[
                H(t) = \int_{0}^{t}\Big[t(E_{n+1}-t)(E_n-s) - s(E_n-t)(E_{n+1}-s)\Big]d\mu_n(s).
            \]
            Expanding the integrand, writing $V_n := \var_{\mu_n}(\tau)$, and noticing that (because $\mu_{n+1}$ is a size-bias of $\mu_n$) $E_n E_{n+1} = V_n + E_n^2$, one finds
            \[
                H(t) = \int_{0}^{t}(t-s)\Big[E_n E_{n+1} - (t+s)E_n + ts\Big]d\mu_n(s) = \int_{0}^{t}(t-s)\Big[V_n + (E_n - t)(E_n - s)\Big]d\mu_n(s).
            \]
            Then it is straightforward to check:
            $t - s > 0$ for $0<s<t$;
            $E_n - t \geq 0$ and $E_n - s > 0$ for $s < t \leq E_n$;
            $V_n > 0$ as a variance of a non-degenerate random variable.
        \end{case}

        \begin{case}{2}{$E_n<t\leq E_{n+1}$}
            Write $H(t)$ simply as
            \[
                H(t) = t(E_{n+1}-t)I_n(t) + E_n(t - E_n)I_{n+1}.
            \]
            We have already shown that $I_n(t) > 0$ for $0<t<1$ and all $n$, and since $E_n < t\leq E_{n+1}$, both terms in the sum are positive.
        \end{case}

        \begin{case}{3}{$E_{n+1}<t\leq 1$}
            Write $H(t)$ under a single integral, but now using the second definition (\ref{eqn:In-def2}):
            \[
                H(t) = \int_{t}^{1}(s-t)\Big[V_n + (E_n - t)(E_n - s)\Big]d\mu_n(s).
            \]
            As before:
            $s - t > 0$ for $t < s < 1$;
            $(E_n - t) < 0$ and $E_n - s < 0$, so their product is positive; 
            $V_n > 0$.
        \end{case} 

        Applying (\ref{eqn:cov-mu_n-change-of-measure-identity}) to (\ref{eqn:A'(n,x)/A(n,x)-as-covariance}) with $f(t) = \partial_x g_x(t)$, since $\partial_t\partial_x g_x(t) = \partial_x\partial_t g_x(t)$, yields
        \[
            \frac{A'(n,x)}{A(n,x)} 
            = \int_{0}^{1}\frac{(\partial_x g_x(t))'}{g'_x(t)}\frac{g'_x(t)I_n(t)}{K_{n,x}}\dt 
            = \e_{\eta_{n,x}}(\Psi(\tau)); \quad \Psi(t) := \partial_x \ln\abs{g'_x(t)}. \qedhere
        \]
    \end{proof}

    \begin{lemma}
        $\Psi'(t) < 0$.
    \end{lemma}

    \begin{proof}
        It will be more convenient to work under the change of variables $u = -\ln(t)$.
        Using (\ref{eqn:g'_x(t)}) and applying this change of variables, for which $\partial_t = -e^u\partial_u$, we can express $\Psi(t)$ as:
        \[
            \Psi(u) = u - \frac{N(u)}{D(u)}; \quad N(u) = e^{-Mu}-1+Mue^{-xu}; \quad D(u) = x(1-e^{-Mu})-M(1-e^{-xu}). 
        \]
        Then, writing $\dot{\Psi}$ for differentiation with respect to $u$,
        \[
            \Psi'(t) = -e^u \dot{\Psi}(u).
        \]
        The result will follow from showing that $\dot{\Psi}(u) > 0$, or ${D(u)}^2 > \dot{N}(u)D(u) - N(u)\dot{D}(u)$.
        To begin, the derivatives of $N$ and $D$ are 
        \begin{align}
            \dot{N}(u) &= M\big[(1-e^{-Mu}) - (1-e^{-xu})\big] - Mxue^{-xu}; \\
            \dot{D}(u) &= Mx\big[(1-e^{-xu})-(1-e^{-Mu})\big].
        \end{align}
        Next, we write ${D(u)}^2 - \dot{N}(u)D(u) + N(u)\dot{D}(u)$ explicitly and collect terms:
        \begin{align*}
            & x^2{(1-e^{-Mu})}^2 + (M^2 - 2Mx)(1-e^{-Mu})(1-e^{-xu}) + Mx(x-M)ue^{-xu}(1-e^{-Mu}) \\
            & \qquad = e^{-xu}(1-e^{-Mu})\Big[x^2(e^{xu}-e^{(x-M)u})+(M^2 - 2Mx)(e^{xu}-1)+Mx(x-M)u\Big] \\
            & \qquad = e^{-xu}(1-e^{-Mu})\Big[{(x-M)}^2 e^{xu} - x^2e^{(x-M)u} + (2Mx-M^2) + Mx(x-M)u\Big].
        \end{align*}
        The leading terms are clearly positive. 
        Write 
        \[
            H(u) := {(x-M)}^2 e^{xu} - x^2e^{(x-M)u} + (2Mx-M^2) + Mx(x-M)u
        \]
        and note
        \[
            \dot{H}(u) = x(x-M)\Big[(x-M)e^{xu} - x e^{(x-M)u} + M\Big], \quad \ddot{H}(u) = x^2{(x-M)}^2 e^{xu}(1-e^{-Mu})
        \]
        and that $H(0) = \dot{H}(0) = 0$.
        Clearly, $\ddot{H}(u) > 0$ and so $\dot{H}(u)$ and, in turn, $H(u)$ are positive as well.
        It follows that ${D(u)}^2 > \dot{N}(u)D(u) - N(u)\dot{D}(u)$ and therefore $\dot{\Psi}(u) > 0$ as desired.
    \end{proof}

    As stated, the result follows from $\eta_{n+1,x}\geq\st\eta_{n,x}$ and $\Psi'(t)<0$, for $x\neq M$.

    To show that $R(n,x)$ has a continuous extension to $x = M$, write $x = M + \delta$.
    Taking the Taylor series of $\psi(n-x)=\psi(n-M-\delta)$ about $n-M$, 
    \begin{align*}
        Q(n,x) &= \frac{\psi(n) - [\psi(n-M) - \delta\psi'(n-M) + \bigO{\delta^2}]}{\psi(n) - \psi(n-M)} \\
        &= 1 + \delta\cdot\frac{\psi'(n-M)}{\psi(n)-\psi(n-M)} + \bigO{\delta^2} \\
        &= 1 + \delta C(n) + \bigO{\delta^2},
    \end{align*}
    where $C(n) := \psi'(n-M)/[\psi(n)-\psi(n-M)]$.
    In turn, writing $\Delta C(n) = C(n+1)-C(n)$,
    \[
        A(n,x) = Q(n+1,x)-Q(n,x) = \delta\Delta C(n) + \bigO{\delta^2},
    \]
    and 
    \[
        R(n,x) = \frac{A(n,x)}{A(n+1,x)} = \frac{\delta\Delta C(n) + \bigO{\delta^2}}{\delta\Delta C(n+1)+\bigO{\delta^2}} = \frac{\Delta C(n) + \bigO{\delta}}{\Delta C(n+1) + \bigO{\delta}} \rightarrow \frac{\Delta C(n)}{\Delta C(n+1)}
    \]
    as $\delta\rightarrow 0$.
\end{proof}

\begin{proof}[Proof of Theorem~\ref{thm:recursion-asymptotics}\ref{item:recursion-asymptotics-error-bound}]
    Note that $s(n,i) = H_{n-1}-H_{n-i-1}$.
    From (\ref{eqn:xn-recursion}), dropping the normalization factor $1/H_{n-1}$, we have 
    \begin{equation}\label{eqn:balance-formula}
        \sum_{i=1}^{n-1}(H_{n-1}-H_{n-i-1})(x_{i+1}-x_i) = 0.
    \end{equation}
    Recall that $x_i = 0$ for $1\leq i\leq M$, except for $x_k > 0$, so for $k < M$, $x_{M+1}-x_M>0$.

    Lemma~\ref{lem:monotonicity} showed that either the tail of $x_n$ is monotone increasing for $n\geq M+1$ or the sequence of differences $x_{n+1}-x_n$ is initially positive, changes sign exactly once, and the tail of $x_n$ is strictly monotone decreasing.    
    For $1\leq k < M$, we derive the error bound separately for each of these cases.

    \begin{case}{1}{Tail of $x_n$ strictly monotone decreasing}
        Suppose that $x_{i+1}-x_i \geq 0$ for $i=M+1, \dots, N$ and $<0$ otherwise, with $M\leq N < \infty$.
        For $k,M\geq 2$ and $k < M$, the terms of (\ref{eqn:balance-formula}) up to $i=M-1$ are zero except for
        \[
            (H_{n-1}-H_{n-k-1})(x_{k+1}-x_k) + (H_{n-1}-H_{n-k})(x_k - x_{k-1}) = x_k(H_{n-k-1} - H_{n-k}) = \frac{-x_k}{n-k}.
        \]
        Then, decomposing the sum as
        \[
            \frac{-x_k}{n-k} + \sum_{i=M}^{N}(H_{n-1}-H_{n-i-1})(x_{i+1}-x_i) + \sum_{i=N+1}^{n-1}(H_{n-1}-H_{n-i-1})(x_{i+1}-x_i)
        \]
        we notice that the first term is negative, the first sum over $i=M,\dots,N$ non-negative, and the second sum negative.
        Recall that $x_n \leq x_k$ for all $n \geq M+1$.
        Then we have, below separating the $x_{N+1}$ term from the sum and grouping the remaining telescoping terms,
        \begin{align*}
            \sum_{i=N+1}^{n-1}(H_{n-1}-H_{n-i-1})(x_i-x_{i+1}) &= -\frac{x_k}{n-k}  + \sum_{i=M}^{N}(H_{n-1}-H_{n-i-1})(x_{i+1}-x_i) \\
            &= (H_{n-1}-H_{n-N-1})x_{N+1} - \frac{x_k}{n-k} - \sum_{i=M+1}^{N}\frac{x_i}{n-i} \\
            &\leq \frac{N}{n-N}x_k
        \end{align*}
        where in the final line we have simply dropped the negative terms in the upper bound and used that $H_{n-1}-H_{n-N-1}\leq \frac{N}{n-N}$.
        We can now bound the error in the same manner as~\cite{Brandenberger_Chin_Mossel_2026}:
        since $H_{n-1}-H_{n-i-1}\geq H_{n-1}-H_{\lfloor n/2\rfloor-1}\geq\ln(2)$ for $\lceil n/2\rceil\leq i\leq n-1$ we have whenever $\lceil n/2\rceil > N$ 
        \begin{align*}
            \abs{x_{\lceil n/2\rceil} - x_n} &\leq \frac{1}{\ln(2)}\sum_{i=\lceil n/2\rceil}(H_{n-1}-H_{n-i-1})\abs{x_i-x_{i+1}} \\
            &= \frac{1}{\ln(2)}\sum_{i=\lceil n/2\rceil}(H_{n-1}-H_{n-i-1})(x_i-x_{i+1}) \\
            &\leq \frac{1}{\ln(2)}\frac{N}{n-N}x_k
        \end{align*}
        since $x_i - x_{i+1} > 0$ for $i\geq \lceil n/2\rceil > N$.
        Writing $x_n - x = (x_n - x_{2n}) + (x_{2n} - x_{4n}) + \dots$, we find 
        \[
            \abs{x_n - x} \leq \frac{N x_k}{\ln(2)}\sum_{j=0}^{\infty}\frac{1}{2^{j}n-N} \leq \frac{1}{\ln(2)}\frac{N x_k}{n-N}\sum_{j=0}^{\infty}\frac{1}{2^j} = \frac{2}{\ln(2)}\frac{N}{n-N}x_k
        \]
        for $\lceil n/2\rceil>N$.
        Lastly, for $n\geq 2N$, $n-N\geq n/2 \geq (n-M)/2$, so $N/(n-N) \leq 2N/(n-M)$, and for $M < n < 2N$, we have simply $x_n - x \leq x_n \leq x_k$ while $2 N x_k / (n - M) \geq 2 N x_k / n > x_k$.
        In either case we obtain 
        \[
            \abs{x_n - x} \leq \frac{4}{\ln(2)}\frac{N}{n-M}x_k
        \]
        for $2\leq k < M$ and some finite $N\geq M$.
        (Note that we are weakening our upper here in order to avoid having $n-N = 0$ for $n = N$.)
        When instead $k = 1$, we decompose (\ref{eqn:balance-formula}) as 
        \[
            -x_k(H_{n-1}-H_{n-k-1}) + \sum_{i=M}^{N}(H_{n-1}-H_{n-i-1})(x_{i+1}-x_i) + \sum_{i=N+1}^{n-1}(H_{n-1}-H_{n-i-1})(x_{i+1}-x_i)
        \]
        and because the first term, still negative, was simply dropped in the previous argument, we obtain the same upper bound when $1 = k < M$.
    \end{case}

    \begin{case}{2}{Tail of $x_n$ monotone increasing}
        Now consider when the sequence is monotone increasing, that is, when $N$ is infinite.
        Then for $2\leq k < M$ we instead decompose (\ref{eqn:balance-formula}) as 
        \[
            \frac{-x_k}{n-k} + \sum_{i=M}^{n-1}(H_{n-1}-H_{n-i-1})(x_{i+1}-x_i).
        \]
        The sum has non-negative terms, and we find simply
        \[
            \sum_{i=M}^{n-1}(H_{n-1}-H_{n-i-1})(x_{i+1}-x_i) = \frac{x_k}{n-k} \leq \frac{k}{n-k}x_k.
        \]
        Using a similar bound to before on $\abs{x_{\lceil n/2\rceil} - x_n}$ whenever $\lceil n/2\rceil \geq M$, that is, 
        \[
            \abs{x_{\lceil n/2\rceil} - x_n} \leq \frac{1}{\ln(2)}\frac{k}{n-k}x_k,
        \] 
        we find
        \[
            \abs{x_n - x} \leq \frac{x_k}{\ln(2)}\sum_{j=0}^{\infty}\frac{1}{2^{j}n-k} \leq \frac{1}{\ln(2)}\frac{x_k}{n-k}\sum_{j=0}^{\infty}\frac{1}{2^j} \leq \frac{2}{\ln(2)}\frac{k}{n-k}x_k
        \]
        for $\lceil n/2\rceil\geq M$.
        For $M\leq n< 2M$, we again have the trivial bound $x_n - x \leq x_n \leq x_k \leq 2 M x_k / (n - M)$.
        As before, we loosen the bound to 
        \[
            \abs{x_n - x} \leq \frac{4}{\ln(2)}\frac{M}{n-M}x_k,
        \]
        for $2\leq k < M$, to match our previous bound.
        When $1=k<M$, we instead write (\ref{eqn:balance-formula}) as 
        \[
            -x_k(H_{n-1}-H_{n-k-1}) + \sum_{i=M}^{n-1}(H_{n-1}-H_{n-i-1})(x_{i+1}-x_i),
        \]
        from which
        \[
            \sum_{i=M}^{n-1}(H_{n-1}-H_{n-i-1})(x_{i+1}-x_i) = x_k(H_{n-1}-H_{n-k-1}) \leq \frac{k}{n-k}x_k,
        \]
        and the same upper bound follows.
    \end{case}

    When $1 = k = M$, we found in Lemma~\ref{lem:monotonicity} that the sequence is constant, so the error bound holds trivially.
\end{proof}

\begin{proof}[Proof of Corollary~\ref{cor:recursion-asymptotics}]
    Let $(x_n^{(k)})$ be the sequence determined by (\ref{eqn:xn-recursion}) with the initial conditions $x_n^{(k)}=0$ for $n\leq M$ except $x_k^{(k)} = x_k > 0$.
    Then as before 
    \[
        x_{M+1}^{(k)} = p(M+1,k)x_k
    \]
    and by linearity of (\ref{eqn:xn-recursion-cbst}) we find 
    \[
        x_{M+1} = \sum_{k=1}^{M}p(M+1,k)x_k = \sum_{k=1}^{M}x_{M+1}^{(k)}.
    \]
    It follows by induction (and again linearity) that $x_n = \sum_{k=1}^{M}x_n^{(k)}$ for all $n\geq M+1$.
    Define $x^{(k)} = \lim_{n\rightarrow\infty}x_n^{(k)}$ and $x := \sum_{k=1}^{M}x^{(k)}$, and let $N_k$ be defined as $N$ before for the sequence $x_n^{(k)}$.
    Then by Theorem~\ref{thm:recursion-asymptotics}\ref{item:recursion-asymptotics-error-bound} we have
    \[
        \abs{x_n - x} \leq \sum_{k=1}^{M}\big|x_n^{(k)} - x^{(k)}\big| \leq \frac{4}{\ln(2)}\frac{N'}{n-M}\sum_{k=1}^{M}x_k
    \]
    where $N' := \max_{1\leq k\leq M} N_k$.
\end{proof}




\bibliographystyle{amsplain}
\bibliography{refs}

@article{Aldous_1991,
	title={Asymptotic Fringe Distributions for General Families of Random Trees},
	volume={1},
	ISSN={1050-5164, 2168-8737},
	DOI={10.1214/aoap/1177005936},
	number={2},
	journal={The Annals of Applied Probability},
	publisher={Institute of Mathematical Statistics},
	author={Aldous, David},
	year={1991},
	month=may,
	pages={228-266}
}

@inbook{Aldous_1996,
	address={New York, NY},
	series={The IMA Volumes in Mathematics and its Applications},
	title={Probability Distributions on Cladograms},
	volume={76},
	ISBN={978-1-4612-6881-9},
	url={http://link.springer.com/10.1007/978-1-4612-0719-1_1},
	DOI={10.1007/978-1-4612-0719-1_1},
	booktitle={Random Discrete Structures},
	publisher={Springer New York},
	author={Aldous, David},
	year={1996},
	pages={1-18},
	collection={The IMA Volumes in Mathematics and its Applications},
	language={en}
}

@article{Aldous_Janson_2025,
	title={The critical beta-splitting random tree {IV}: Mellin analysis of leaf height},
	volume={30},
	ISSN={1083-6489, 1083-6489},
	DOI={10.1214/25-EJP1332},
	number={none},
	journal={Electronic Journal of Probability},
	publisher={Institute of Mathematical Statistics and Bernoulli Society},
	author={Aldous, David and Janson, Svante},
	year={2025},
	month=jan,
	pages={1-39},
	language={en}
}

@article{Aldous_Pittel_2025,
	title={The critical beta-splitting random tree {I}: Heights and related results},
	volume={35},
	ISSN={1050-5164, 2168-8737},
	DOI={10.1214/24-AAP2112},
	number={1},
	journal={The Annals of Applied Probability},
	publisher={Institute of Mathematical Statistics},
	author={Aldous, David and Pittel, Boris},
	year={2025},
	month=feb,
	pages={158-195},
	language={en}
}

@misc{Aldous_Janson_2024a,
	title={The Critical Beta-splitting Random Tree {II}: Overview and Open Problems},
	url={http://arxiv.org/abs/2303.02529},
	DOI={10.48550/arXiv.2303.02529},
	note={arXiv:2303.02529},
	publisher={arXiv},
	author={Aldous, David J. and Janson, Svante},
	year={2024},
	month=jul
}

@misc{Aldous_Janson_2024b,
	title={The Critical Beta-splitting Random Tree {III}: The exchangeable partition representation and the fringe tree},
	url={http://arxiv.org/abs/2412.09655},
	DOI={10.48550/arXiv.2412.09655},
	note={arXiv:2412.09655},
	publisher={arXiv},
	author={Aldous, David J. and Janson, Svante},
	year={2024},
	month=dec
}

@misc{Aldous_Janson_Li_2024,
	title={The Harmonic Descent Chain},
	url={http://arxiv.org/abs/2405.05102},
	DOI={10.48550/arXiv.2405.05102},
	note={arXiv:2405.05102},
	publisher={arXiv},
	author={Aldous, David J. and Janson, Svante and Li, Xiaodan},
	year={2024},
	month=may
}

@misc{Brandenberger_Chin_Mossel_2026,
	title={Asymptotics for the harmonic descent chain and applications to critical beta-splitting trees},
	url={http://arxiv.org/abs/2505.24821},
	DOI={10.48550/arXiv.2505.24821},
	note={arXiv:2505.24821},
	publisher={arXiv},
	author={Brandenberger, Anna and Chin, Byron and Mossel, Elchanan},
	year={2026},
	month=sep
}

@book{Chen_Goldstein_Shao_2011,
	address={Berlin, Heidelberg},
	series={Probability and Its Applications},
	title={Normal Approximation by {S}tein's Method},
	ISBN={978-3-642-15006-7},
	DOI={10.1007/978-3-642-15007-4},
	publisher={Springer-Verlag Berlin Heidelberg},
	author={Chen, Louis H. Y. and Goldstein, Larry Joel and Shao, Qi-Man},
	year={2011},
	collection={Probability and Its Applications},
	language={eng}
}

@article{Erdélyi_Tricomi_1951,
	title={The asymptotic expansion of a ratio of gamma functions.},
	volume={1},
	ISSN={0030-8730},
	number={1},
	journal={Pacific Journal of Mathematics},
	publisher={Pacific Journal of Mathematics, A Non-profit Corporation},
	author={Erdélyi, A. and Tricomi, F. G.},
	year={1951},
	month=jan,
	pages={133-142},
	language={en}
}

@article{Flajolet_Sedgewick_1995,
	title={Mellin transforms and asymptotics: Finite differences and {R}ice's integrals},
	volume={144},
	ISSN={0304-3975},
	DOI={10.1016/0304-3975(94)00281-M},
	number={1},
	journal={Theoretical Computer Science},
	author={Flajolet, Philippe and Sedgewick, Robert},
	year={1995},
	month=jun,
	pages={101-124}
}

@article{Goldstein_Reinert_1997,
	title={Stein's method and the zero bias transformation with application to simple random sampling},
	volume={7},
	ISSN={1050-5164, 2168-8737},
	DOI={10.1214/aoap/1043862419},
	number={4},
	journal={The Annals of Applied Probability},
	publisher={Institute of Mathematical Statistics},
	author={Goldstein, Larry and Reinert, Gesine},
	year={1997},
	month=nov,
	pages={935-952},
	language={en}
}

@article{Iksanov_2025,
	title={The harmonic descent chain and regenerative composition structures},
	volume={30},
	ISSN={1083-589X, 1083-589X},
	DOI={10.1214/25-ECP661},
	number={none},
	journal={Electronic Communications in Probability},
	publisher={Institute of Mathematical Statistics and Bernoulli Society},
	author={Iksanov, Alexander},
	year={2025},
	month=jan,
	pages={1-3},
	language={en}
}

@article{Kolesnik_2025,
	title={Critical beta-splitting, via contraction},
	volume={30},
	ISSN={1083-589X, 1083-589X},
	DOI={10.1214/25-ECP658},
	number={none},
	journal={Electronic Communications in Probability},
	publisher={Institute of Mathematical Statistics and Bernoulli Society},
	author={Kolesnik, Brett},
	year={2025},
	month=jan,
	pages={1-14},
	language={en}
}

@inbook{Murty_2023,
	address={Cham},
	title={Irrationality of Zeros of the Digamma Function},
	ISBN={978-3-031-31617-3},
	url={https://doi.org/10.1007/978-3-031-31617-3_15},
	DOI={10.1007/978-3-031-31617-3_15},
	booktitle={Number Theory in Memory of Eduard Wirsing},
	publisher={Springer International Publishing},
	author={Murty, M. Ram},
	year={2023},
	pages={237-243},
	language={en}
}

@article{Régnier_1989,
	title={A limiting distribution for {Q}uicksort},
	volume={23},
	ISSN={0988-3754},
	number={3},
	journal={RAIRO - Theoretical Informatics and Applications - Informatique Théorique et Applications},
	publisher={EDP-Sciences},
	author={Régnier, Mireille},
	year={1989},
	pages={335-343},
	language={eng}
}

@article{Rösler_1991,
	title={A limit theorem for "{Q}uicksort"},
	volume={25},
	ISSN={0988-3754, 1290-385X},
	DOI={10.1051/ita/1991250100851},
	number={1},
	journal={RAIRO - Theoretical Informatics and Applications},
	author={Rösler, Uwe},
	year={1991},
	pages={85-100},
	language={en}
}

@book{Shaked_Shanthikumar_2007,
	address={New York},
	series={Springer series in statistics},
	title={Stochastic orders},
	ISBN={978-0-387-32915-4},
	callNumber={QA273.6 .S523 2007},
	publisher={Springer},
	author={Shaked, Moshe and Shanthikumar, J. George},
	year={2007},
	collection={Springer series in statistics}
}

@book{Whittaker_Watson_1927,
	address={Cambridge New York},
	edition={Fourth edition, reprintedition},
	series={Cambridge mathematical library},
	title={A course of modern analysis: an introduction to the general theory of infinite processes and of analytic functions, with an account of the principal transcendental functions},
	ISBN={978-0-511-60875-9},
	publisher={Cambridge University Press},
	author={Whittaker, Edmund T. and Watson, G. N.},
	year={1927},
	collection={Cambridge mathematical library},
	language={eng}
}











\begin{acks}
  The author used Claude Sonnet 5.5 during the preparation of this manuscript for proofreading and copyediting purposes.
  LLMs were not otherwise used to generate the content of this work, and all proofs and writing in the manuscript were produced by the author.
\end{acks}


\end{document}